\documentclass[11pt,reqno]{amsart}
\usepackage{amsmath,amsfonts, amssymb, amsthm}
  \usepackage{paralist}
  \usepackage{graphics} 
  \usepackage{epsfig} 
\usepackage{graphicx}  \usepackage{epstopdf}
 \usepackage[colorlinks=true]{hyperref}
\hypersetup{urlcolor=blue, citecolor=red}

\newtheorem{theorem}{Theorem}

\newtheorem{conjecture}{Conjecture}

\theoremstyle{definition}

\newtheorem{remark}{Remark}

\def\pmod #1{\ ({\rm{mod}}\ #1)}
\def\Z{\Bbb Z}
\def\N{\Bbb N}

\def\l{\left}
\def\r{\right}
\def\bg{\bigg}
\def\({\bg(}
\def\){\bg)}
\def\t{\text}
\def\f{\frac}
\def\mo{{\rm{mod}\ }}
\def\pmod#1{\ (\mo\ #1)}

\def\ls{\leqslant}

\def\sm{\setminus}
\def\bi{\binom}

\def\eq{\equiv}

\def\Proof{\noindent{\it Proof}}
\def\Ack{\medskip\noindent {\bf Acknowledgments}}

\begin{document}
\hbox{Included in: M. B. Nathanson (ed.), Combinatorial and Additive Number Theory VI,}
 \hbox{Springer Proc. Math. Stat. 464, Springer, Cham, 2025, pp. 413--460.}
\medskip

\title[Series with summands involving harmonic numbers]
      {Series with summands \\ involving harmonic numbers}
\author[Zhi-Wei Sun]{Zhi-Wei Sun}


\address{School of Mathematics, Nanjing
University, Nanjing 210093, People's Republic of China}
\email{{\tt zwsun@nju.edu.cn}
\newline\indent
{\it Homepage}: {\tt http://maths.nju.edu.cn/\lower0.5ex\hbox{\~{}}zwsun}}

\keywords{Harmonic numbers, series for $\pi$,
Dirichlet $L$-functions, combinatorial identities, congruences.
\newline \indent 2020 {\it Mathematics Subject Classification}. Primary 11B65, 05A19; Secondary 11A07, 11B68. \newline\indent
The initial version of this paper was posted to arXiv in Oct. 2022 as a preprint with the ID arXiv:2210.07238.}

\begin{abstract}
For each positive integer $m$, the $m$th order harmonic numbers are given by
$$H_n^{(m)}=\sum_{0<k\ls n}\frac1{k^m}\ \ (n=0,1,2,\ldots).$$
We discover exact values of some series involving harmonic numbers of order not exceeding four.
For example, we conjecture that $$\sum_{k=0}^\infty(6k+1)\frac{\binom{2k}k^3}{256^k}\left(H_{2k}^{(3)}-\f{7}{64}H_{k}^{(3)}\right)
 =\frac{25\zeta(3)}{8\pi}-G,$$
 where $G$ denotes the Catalan constant $\sum_{k=0}^\infty(-1)^k/(2k+1)^2$.
This paper contains $70$ conjectures posed by the author during 2022--2023.
\end{abstract}
\maketitle

\section{Introduction}
\setcounter{equation}{0}
 \setcounter{conjecture}{0}
 \setcounter{theorem}{0}
 \setcounter{proposition}{0}

 The usual harmonic numbers are those rational numbers
 $$H_n=\sum_{0<k\ls n}\f1k\quad \ (n=0,1,2,\ldots).$$
 For each $m\in\Z^+=\{1,2,3,\ldots\}$, the harmonic numbers of order $m$ are defined by
 $$H_n^{(m)}=\sum_{0<k\ls n}\f1{k^m}\quad \ (n\in\N=\{0,1,2,\ldots\}).$$
 For any $m,n\in\Z^+$, we clearly have
 $$\sum_{k=1}^n\f1{(2k-1)^m}=\sum_{k=1}^{2n}\f1{k^m}-\sum_{j=1}^n\f1{(2j)^m}=H_{2n}^{(m)}-\f1{2^m}H_n^{(m)}.$$
 J. Wolstenholme \cite{W} established two fundamental congruences for harmonic numbers:
 $$H_{p-1}\eq0\pmod {p^2}\ \ \t{and}\ \ H_{p-1}^{(2)}\eq0\pmod p$$
 for any prime $p>3$.  For series and congruences involving harmonic numbers,
 one may consult \cite{S12,S15e,SZ}, \cite[Section 10.5]{S-Book}, and the recent preprint \cite{Au22} solving various conjectures of the author.

 In 2012, K. N. Boyadzhiev \cite{Boy} proved that
 \begin{equation}\label{CHk}\sum_{k=0}^\infty\bi {2k}kH_kx^k=\f{2}{\sqrt{1-4x}}\log\f{1+\sqrt{1-4x}}{2\sqrt{1-4x}}
 \quad \t{for}\ x\in\l(-\f14,\f14\r).
 \end{equation} In 2016, H. Chen \cite{Chen} deduced that
 \begin{equation}\label{CH2k}\sum_{k=0}^\infty\bi {2k}kH_{2k}x^k=\f{1}{\sqrt{1-4x}}\log\f{1+\sqrt{1-4x}}{2(1-4x)}
 \quad \t{for}\ x\in\l(-\f14,\f14\r).
 \end{equation}

 It is well known that
$$2\arcsin\f x2=\sum_{k=0}^\infty\f{\bi{2k}kx^{2k+1}}{(2k+1)16^k}\quad \t{for}\ |x|\ls2;$$
in particular,
 $$\sum_{k=0}^\infty\f{\bi{2k}k}{(2k+1)16^k}=\f{\pi}3
 \ \ \t{and}\ \ \sum_{k=0}^\infty\f{\bi{2k}k}{(2k+1)8^k}=\f{\pi}{2\sqrt2}.$$
 The author \cite[Theorem 1.1(ii)]{S11} determined
 $$\sum_{k=0}^{(p-3)/2}\f{\bi{2k}k}{(2k+1)16^k}\ \ \t{and}\ \ \sum_{k=(p+1)/2}^{p-1}\f{\bi{2k}k}{(2k+1)16^k}
 $$
 modulo $p^2$ for any prime $p>3$. By \cite{BC}, we have
\begin{equation}\label{cube}\l(\mathrm{arcsin}\f x2\r)^3=3\sum_{k=0}^\infty\f{\bi{2k}kx^{2k+1}}{(2k+1)16^k}\sum_{0\ls j<k}\f1{(2j+1)^2}
\end{equation}
and
\begin{equation}\label{four}\l(\mathrm{arcsin}\f x2\r)^4=\f 32\sum_{k=1}^\infty\f{H_{k-1}^{(2)}x^{2k}}{k^2\bi{2k}k}
\end{equation}
for $|x|\ls 2$.
In particular,
\begin{align*}\sum_{k=0}^\infty\f{\bi{2k}k}{(2k+1)16^k}\l(H_{2k}^{(2)}-\f14H_k^{(2)}\r)&=\f{\pi^3}{648},
\\\sum_{k=0}^\infty\f{\bi{2k}k}{(2k+1)8^k}\l(H_{2k}^{(2)}-\f14H_k^{(2)}\r)&=\f{\sqrt2\,\pi^3}{384}.
\end{align*}
and
$$\sum_{k=1}^\infty\f{H_{k-1}^{(2)}}{k^2\bi{2k}k}=\f{\pi^4}{1944}.$$

 In view of \eqref{cube}, we have the following result.

 \begin{theorem} If $|x|<2$, then
 \begin{equation}\label{two}\f{(\arcsin(x/2))^2}{\sqrt{4-x^2}}=\sum_{k=1}^\infty\f{\bi{2k}kx^{2k}}{16^k}\l(H_{2k}^{(2)}-\f14H_k^{(2)}\r).
 \end{equation}
 \end{theorem}
\Proof. By taking derivatives of both sides of \eqref{cube}, we get
$$3\l(\arcsin\f x2\r)^2\times\f{1/2}{\sqrt{1-(x/2)^2}}=3\sum_{k=0}^\infty\f{\bi{2k}kx^{2k}}{16^k}\sum_{0\ls j<k}\f1{(2j+1)^2}$$
and hence
$$\f{(\arcsin(x/2))^2}{\sqrt{4-x^2}}=\sum_{k=1}^\infty\f{\bi{2k}kx^k}{16^k}\sum_{j=1}^k\f1{(2j-1)^2},$$
which is equivalent to \eqref{two}. \qed

 Motivated by the Ramanujan series
 \begin{equation}\label{-512R}\sum_{k=0}^\infty(6k+1)\f{\bi{2k}k^3}{(-512)^k}=\f{2\sqrt2}{\pi}
 \end{equation}
 (cf. \cite{R}),
 L. Long \cite{Long} conjectured the congruence
 \begin{equation}\label{L-cong}\sum_{k=0}^{(p-1)/2}(6k+1)\f{\bi{2k}k^3}{(-512)^k}\sum_{j=1}^k\l(\f1{(2j-1)^2}-\f1{16j^2}\r)\eq0\pmod p\end{equation}
 for any odd prime $p$, which was confirmed by H. Swisher \cite{Sw} in 2015.
Note that \eqref{L-cong} can be rewritten as
$$\sum_{k=0}^{(p-1)/2}(6k+1)\f{\bi{2k}k^3}{(-512)^k}\l(H_{2k}^{(2)}-\f5{16}H_k^{(2)}\r)\eq0\pmod p.$$

In 2022 C. Wei \cite{WeiRJ} deduced the two identities
$$\sum_{k=0}^{\infty}(6k+1)\f{\bi{2k}k^3}{(-512)^k}\l(H_{2k}^{(2)}-\f5{16}H_k^{(2)}\r)=-\f{\sqrt2}{48}\pi$$
and
$$\sum_{k=0}^{\infty}(6k+1)\f{\bi{2k}k^3}{256^k}\l(H_{2k}^{(2)}-\f5{16}H_k^{(2)}\r)=\f{\pi}{12}$$
conjectured by V.J.W. Guo and X. Lian \cite{GL}, as well as their $q$-analogues.

Motivated by Bauer's series
\begin{equation}\label{Bauer}\sum_{k=0}^\infty(4k+1)\f{\bi{2k}k^3}{(-64)^k}=\f2{\pi}
\end{equation}
and Ramanujan's series
\begin{equation}\label{48^2R}\sum_{k=0}^\infty(8k+1)\f{\bi{2k}k^2\bi{4k}{2k}}{48^{2k}}=\f{2\sqrt3}{\pi},
\end{equation}
Wei and G. Ruan \cite{W22} proved the two new identities
$$\sum_{k=1}^\infty(4k+1)\f{\bi{2k}k^3}{(-64)^k}\sum_{j=1}^{2k}\f{(-1)^j}{j^2}=\f{\pi}{12}$$
and
$$\sum_{k=1}^\infty(8k+1)\f{\bi{2k}k^2\bi{4k}{2k}}{48^{2k}}\sum_{j=1}^k\l(\f1{(2j-1)^2}-\f1{36j^2}\r)
=\f{\sqrt3\,\pi}{54},$$
i.e.,
\begin{equation}\label{64}
\sum_{k=0}^\infty(4k+1)\f{\bi{2k}k^3}{(-64)^k}\l(H_{2k}^{(2)}-\f12H_k^{(2)}\r)=-\f{\pi}{12}
\end{equation}
and
\begin{equation}\label{48^2}\sum_{k=0}^\infty(8k+1)\f{\bi{2k}k^2\bi{4k}{2k}}{48^{2k}}
\l(H_{2k}^{(2)}-\f 5{18}H_k^{(2)}\r)=\f{\sqrt3\,\pi}{54}.
\end{equation}

In 1997 van Hamme \cite{vH} thought that series for powers of $\pi=\Gamma(1/2)^2$ should have their $p$-adic analogues
involving the $p$-adic Gamma function $\Gamma_p(x)$, where $p$ is an odd prime.
Note that for any odd prime $p$ we have
$$\Gamma_p\l(\f12\r)^2=(-1)^{(p+1)/2}=-\l(\f{-1}p\r),$$
where $(\f{\cdot}p)$ denotes the Legendre symbol.
The author \cite{S20,S21,S23} found many new series for powers of $\pi$ motivated by related congruences.
However, van Hamme's philosophy fails
for some Ramanujan-type series for $1/\pi$. For example, T. Huber, D. Schultz and D. Ye \cite{HSY}
used modular forms to obtain that
$$\sum_{k=0}^\infty(6k+1)\f{a_k}{16^k}=\f{16}{\pi},$$
where $a_0=1,\ a_1=4,\ a_2=20$ and
$$(n+1)^3a_{n+1}=4(2n+1)(2n^2+2n+1)a_n-16n(4n^2+1)a_{n-1}+8(2n-1)^3a_{n-2}$$
for all $n=2,3,\ldots$; but for a general odd prime $p$ we even cannot find any pattern for
$\sum_{k=0}^{p-1}(6k+1)a_k/16^k$ modulo $p$.

The Bernoulli numbers $B_0,B_1,B_2,\ldots$ are defined by
 $$\f x{e^x-1}=\sum_{n=0}^\infty B_n\f{x^n}{n!}\ \ \l(0<|x|<2\pi\r).$$
 Equivalently,
 $$B_0=1,\ \t{and}\ \sum_{k=0}^{n}\bi{n+1}kB_k=0\ \ \t{for}\ n=1,2,3,\ldots.$$
 In 1900 J.W.L. Glaiser \cite{Gl} proved that
 $$H_{p-1}\eq-\f{p^2}3B_{p-3}\pmod{p^3}\ \ \t{and}\ \ H_{p-1}^{(2)}\eq\f23 pB_{p-3}\pmod{p^2}$$
 for any prime $p>3$.
 The Bernoulli polynomials are given by
 $$B_n(x)=\sum_{k=0}^n\bi nk B_kx^{n-k}\ \ (n\in\N).$$
The Euler numbers $E_0,E_1,E_2,\ldots$ are defined by
 $$\f2{e^x+e^{-x}}=\sum_{n=0}^\infty E_n\f{x^n}{n!}\ \ \l(|x|<\f{\pi}2\r).$$
 Clearly $E_{2n+1}=0$ for all $n\in\N$. It is also known that
 $$\sum_{k=0}^n\bi{2n}{2k}E_{2k}=0\ \ \ \t{for each}\ n\in\N.$$
 The Euler polynomials are given by
 $$E_n(x)=\sum_{k=0}^n\bi nk\f{E_k}{2^k}\l(x-\f12\r)^{n-k}\ \ (n\in\N).$$
 The author \cite{S11c,S13} first observed that many Ramanujan-type series
 have corresponding congruences involving Bernoulli or Euler polynomials.

 Now we introduce some notations throughout this paper. The Riemann zeta function is defined by
 $$\zeta(s)=\sum_{n=1}^\infty\f1{n^s}\ \ \ \t{with}\ \Re(s)>1.$$
 The Dirichlet beta function is given by
 $$\beta(m)=\sum_{k=0}^\infty\f{(-1)^k}{(2k+1)^m}\ \ \ \ (m=1,2,3,\ldots).$$
 Note that $G=\beta(2)$ is the Catalan constant. We also adopt the notation
$$K:=L\l(2,\l(\f{-3}{\cdot}\r)\r)=\sum_{n=1}^\infty\f{\l(\f k3\r)}{k^2}$$
with $(\f{-3}{\cdot})$ the Kronecker symbol. For a prime $p$ and an integer $a\not\eq0\pmod p$,  we use $q_p(a)$ to denote the Fermat quotient
 $(a^{p-1}-1)/p$. Many congruences in later sections involve Fermat quotients.

 In Sections 2--4, we will propose 70 new conjectures on series and related congruences
 with summands involving not only harmonic numbers of order at most four, but also
 products of several binomial coefficients.
 All the conjectures have been checked via {\tt Mathematica}.

\section{Series with summands containing \\ one or two binomial coefficients}
 \setcounter{theorem}{0}
 \setcounter{proposition}{0}

 \begin{conjecture}[{\rm 2022-10-12}] {\rm (i)} We have
 \begin{equation}\label{451}\sum_{k=1}^\infty\f{(-1)^{k-1}}{k^3\bi{2k}k}\l(H_{2k-1}^{(2)}-\f{123}{16}H_{k-1}^{(2)}\r)
 =\f{451}{40}\zeta(5)-\f{14}{15}\pi^2\zeta(3).
 \end{equation}

 {\rm (ii)} For any prime $p>5$, we have
 \begin{equation}\sum_{k=1}^{(p-1)/2}\f{(-1)^k}{k^3\bi{2k}k}\l(16H_{2k-1}^{(2)}-123H_{k-1}^{(2)}\r)
 \eq-542B_{p-5}\pmod p
 \end{equation}
 and
 \begin{equation}p\sum_{k=1}^{p-1}\f{(-1)^k}{k^2}\bi{2k}k\l(16H_{2k}^{(2)}-123H_k^{(2)}\r)\eq 192\f{H_{p-1}}{p^2}\pmod{p^2}.
 \end{equation}
 \end{conjecture}
 \begin{remark} In 1979 R. Ap\'ery \cite{Ap} proved the irrationality of $\zeta(3)=\sum_{n=1}^\infty1/n^3$ via the identity
 $$\sum_{k=1}^\infty\f{(-1)^{k-1}}{k^3\bi{2k}k}=\f25\zeta(3).$$
 In 2014 the author \cite{S14JNT} proved the congruence
 $$\sum_{k=1}^{(p-1)/2}\f{(-1)^{k-1}}{k^3\bi{2k}k}\eq 2B_{p-3}\pmod p$$
 for any prime $p>5$.
 The author's conjectural identity (cf. \cite{S15})
 $$\sum_{k=1}^\infty\f{(-1)^{k-1}}{k^3\bi{2k}k}(H_{2k}+4H_k)=\f{2\pi^4}{75}$$
 was proved by W. Chu \cite{Chu} as well as K. C. Au \cite[Prop. 7.14]{Au}.
 After seeing an earlier arXiv version of this paper,  Au \cite[Corollary 2.9]{Au22}
 confirmed the author's conjectural identity \eqref{451}.
 \end{remark}

\begin{conjecture}[{\rm 2022-11-14}] We have the identity
\begin{equation}\label{8^k}\sum_{k=0}^\infty\f{\bi{2k}k}{8^k}\l(H_{2k}^{(3)}-\f18 H_k^{(3)}\r)
=\f{35\sqrt2}{64}\zeta(3)-\f{\sqrt2}8\pi G.
\end{equation}
\end{conjecture}
\begin{remark} Applying \eqref{CHk} and \eqref{CH2k} with $x=1/8$, we see that
$$\sum_{k=0}^\infty\f{\bi{2k}k}{8^k}H_k=-\sqrt2\log(12-8\sqrt2)
\ \t{and}\ \sum_{k=0}^\infty\f{\bi{2k}k}{8^k}H_{2k}=\f{\log(3/2+\sqrt2)}{\sqrt2}.$$
In contrast with \eqref{8^k}, we have
$$\sum_{k=0}^\infty\f{\bi{2k}k}{8^k}\l(H_{2k}^{(2)}-\f14H_k^{(2)}\r)=\f{\pi^2}{16\sqrt2}$$
by applying \eqref{two} with $x=\sqrt2$.
\end{remark}

\begin{conjecture}[{\rm 2022-11-14}] We have the identity
\begin{equation}\label{16k}\sum_{k=0}^\infty\f{\bi{2k}k}{16^k}\l(H_{2k}^{(3)}-\f18 H_k^{(3)}\r)
=\f{2\zeta(3)}{3\sqrt3}-\f{\pi K}8.
\end{equation}
\end{conjecture}
\begin{remark}  Applying \eqref{CHk} and \eqref{CH2k} with $x=1/8$, we see that
$$\sum_{k=0}^\infty\f{\bi{2k}k}{16^k}H_k=-\f2{\sqrt3}\log(84-48\sqrt3)
\ \t{and}\ \sum_{k=0}^\infty\f{\bi{2k}k}{16^k}H_{2k}=\f{\log((7+4\sqrt3)/9)}{\sqrt3}.$$
In contrast with \eqref{16k}, we have
$$\sum_{k=0}^\infty\f{\bi{2k}k}{16^k}\l(H_{2k}^{(2)}-\f14H_k^{(2)}\r)=\f{\pi^2}{36\sqrt3}$$
by applying \eqref{two} with $x=1$.
\end{remark}

\begin{conjecture}[{\rm 2023-05-28}] {\rm (i)} We have
 \begin{equation}\sum_{k=1}^\infty\f{H_{3k}-H_k}{k2^k\bi{3k}k}=\f25(G+\log^22)-\f{\pi^2}{24}
 \end{equation}
 and
 \begin{equation}\sum_{k=1}^\infty\f{H_{3k}-H_k}{k^22^k\bi{3k}k}=\f{11}4\zeta(3)-\f{\pi^2}{24}\log2-\pi G.
 \end{equation}

 {\rm (ii)} For any prime $p>5$ with $p\eq 1\pmod4$, we have
 \begin{equation} p\sum_{k=1}^{p-1}\f{H_{3k}-H_k}{k2^k\bi{3k}k}\eq\f 7{10}q_p(2)\pmod p.
 \end{equation}
 Also, for each odd prime $p$ we have
 \begin{equation}p^2\sum_{k=1}^{p-1}\f{H_{3k}-H_k}{k^22^k\bi{3k}k}\eq-\f{q_p(2)}4\pmod p.
 \end{equation}
 \end{conjecture}
 \begin{remark} The author's conjectural identities
 $$\sum_{k=1}^\infty\f{H_{2k}-H_k}{k2^k\bi{3k}k}=\f3{10}\log^22+\f{\pi}{20}\log2-\f{\pi^2}{60}$$
 and
 $$\sum_{k=1}^\infty\f{H_{2k}-H_k}{k^22^k\bi{3k}k}=\f{33}{32}\zeta(3)+\f{\pi^2}{24}\log2-\f{\pi G}2$$
 (cf. \cite[Conjecture 10.61]{S-Book}) were confirmed by Au \cite{Au} in 2022.
 \end{remark}

 \begin{conjecture}[{\rm 2023-05-28}] {\rm (i)} We have
 \begin{equation}\sum_{k=1}^\infty\f{25k-3}{2^k\bi{3k}k}(H_{3k}-8H_{2k}+7H_k)=2G-2(\pi+9)\log2.
 \end{equation}

 {\rm (ii)} For any odd prime $p$, we have the congruence
 \begin{equation}p^2\sum_{k=1}^{p-1}\f{25k-3}{2^k\bi{3k}k}(H_{3k}-8H_{2k}+7H_k)\eq-\l(\f{-1}p\r)\f 94\pmod p.
 \end{equation}
 \end{conjecture}
 \begin{remark} In 1974 R. W. Gosper announced the identity
 $$\sum_{k=0}^\infty\f{25k-3}{2^k\bi{3k}k}=\f{\pi}2,$$
 an elegant proof of which can be found in \cite{AKP}.
 \end{remark}

 \begin{conjecture}[{\rm 2023-05-28}] {\rm (i)} We have
 \begin{equation}\sum_{k=0}^\infty\f{2^k\bi{3k}k}{27^k}H_{2k}=\f32(1+\sqrt3)\log\l(1+\f{\sqrt3}3\r)-\sqrt3\log2
 \end{equation}
 and
 \begin{equation}\sum_{k=0}^\infty\f{2^k\bi{3k}k}{27^k}H_{3k}=\f{1+\sqrt3}2\l(2\log(1+\sqrt3)-\f{\log3}2\r)-\sqrt3\log2.
 \end{equation}

 {\rm (ii)} For any prime $p>3$, we have
 \begin{equation}\sum_{k=(p+1)/2}^{p-1}\bi{3k}k\l(\f2{27}\r)^k\eq\f{1-(\f 3p)}3\pmod p.
 \end{equation}
 \end{conjecture}
 \begin{remark} For any positive integer $n$, we clearly have
 $$H_n=\sum_{k=0}^{n-1}\f1{k+1}=\sum_{k=0}^{n-1}\int_0^1 t^kdt=\int_0^1\sum_{k=0}^{n-1}t^kdt=\int_0^1\f{1-t^n}{1-t}dt.$$
 Using this trick we can deduce that
 $$\sum_{k=0}^\infty\f{2^k\bi{3k}k}{27^k}H_k=\f34\l((1-\sqrt3)\log4-(1+\sqrt3)\log3\r)+2\sqrt3\log(1+\sqrt3).$$
 \end{remark}

 \begin{conjecture}[{\rm 2023-05-28}] We have
 \begin{equation}\sum_{k=0}^\infty\bi{3k}k\l(\f{3+\sqrt5}{54}\r)^k(H_{3k}-H_{2k})=\phi(\log3-2\log\phi),
 \end{equation}
 where $\phi$ denotes the golden ratio $(1+\sqrt5)/2\approx 1.618\ldots$.
 \end{conjecture}
 \begin{remark} {\tt Mathematica} yields that
 $$\sum_{k=0}^\infty\bi{3k}k\l(\f{4x}{27}\r)^k=\f1{\sqrt{1-x}}\cos\f{\arcsin\sqrt x}3$$
 for any $x\in(-1,1)$. Applying this with $x=((1+\sqrt5)/4)^2$ we obtain that
 $$\sum_{k=0}^\infty\bi{3k}k\l(\f{3+\sqrt5}{54}\r)^k=\f{\cos (\pi/10)}{\sqrt{(5-\sqrt5)/8}}=\f{\sqrt{(5+\sqrt5)/8}}{\sqrt{(5-\sqrt5)/8}}=\phi.$$
 \end{remark}

 \begin{conjecture}[{\rm 2023-05-30}] If $(1-\sqrt2)/2\ls x<1/2$, then
 \begin{equation}\sum_{k=0}^\infty\bi{4k}{2k}\l(\f{x(1-x)}4\r)^k(2H_{4k}-3H_{2k}+H_k)=\f{\sqrt{1-x}}{2x-1}\log(1-x).
 \end{equation}
 \end{conjecture}
 \begin{remark} For any $x\in(-1,1)$, we have
 $$\sum_{k=0}^\infty\bi{4k}{2k}\left(\f x{16}\r)^k=\sqrt{\f{1+\sqrt{1-x}}{2(1-x)}}$$
 which can be proved directly or via {\tt Mathematica}. In particular,
 $$\sum_{k=0}^\infty\bi{4k}{2k}\left(\f3{64}\r)^k=\sqrt3.$$
 For $x\in(-1,1)$, we obviously have
 $$\sum_{k=0}^\infty \bi{2k}kH_k\l(\l(\f x4\r)^k+\l(-\f x4\r)^k\r)=2\sum_{k=0}^\infty\bi{4k}{2k}H_{2k}\l(\f{x}4\r)^{2k}$$
 and
 $$\sum_{k=0}^\infty\bi{2k}k H_{2k}\l(\l(\f x4\r)^k+\l(-\f x4\r)^k\r)=2\sum_{k=0}^\infty\bi{4k}{2k}H_{4k}\l(\f{x}4\r)^{2k},$$
 and hence we may find closed formulas for the two series
 $$\sum_{k=0}^\infty \bi{4k}{2k}H_{2k}\l(\f x4\r)^{2k}
 \ \ \t{and}\ \ \sum_{k=0}^\infty\bi{2k}k H_{4k}\l(\f x4\r)^{2k}$$
 by using \eqref{CHk} and \eqref{CH2k}.
 In particular,
 $$\sum_{k=0}^\infty\bi{4k}{2k}\l(\f{3}{64}\r)^kH_{2k}=2(1+\sqrt3)\log(1+\sqrt3)-(1+3\sqrt3)\log2+\f{\sqrt3-1}2\log3$$
 and
 $$\sum_{k=0}^\infty\bi{4k}{2k}\l(\f{3}{64}\r)^kH_{4k}=(2+\sqrt3)\log(1+\sqrt3)-(1+\sqrt3)\log2+\f{\sqrt3-1}4\log 3.$$
 With the aid of {\tt Mathematica}, we obtain that
 \begin{align*}\sum_{k=0}^\infty\bi{4k}{2k}\left(\f3{64}\r)^kH_k&=\sum_{k=1}^\infty\bi{4k}{2k}\left(\f3{64}\r)^k\int_0^1\f{1-t^k}{1-t}dt
 \\&=\int_0^1\f1{1-t}\sum_{k=1}^\infty\bi{4k}{2k}\left(\f3{64}\r)^k(1-t^k)dt
 \\&=\int_0^1\f1{1-t}\left(\sqrt3-\f{\sqrt{2+\sqrt{4-3t}}}{\sqrt{4-3t}}\r)dt
 \\&=\log\f{2+\sqrt3}3+\sqrt3\log\f{7+4\sqrt3}8
 \\&=(2+4\sqrt3)\log(1+\sqrt3)-\log3-(1+5\sqrt3)\log2.
 \end{align*}
 \end{remark}

 \begin{conjecture} [{\rm 2022-11-14}] We have the identity
 \begin{equation}\label{32^k}\sum_{k=0}^\infty\f{\bi{2k}k^2}{32^k}\l(H_{2k}^{(2)}-\f14H_k^{(2)}\r)
 =\Gamma\l(\f14\r)^2\f{\pi^2-8G}{32\pi\sqrt{\pi}},
 \end{equation}
 where $\Gamma(x)$ is the well-known Gamma function.
 \end{conjecture}
 \begin{remark} In contrast with \eqref{32^k}, {\tt Mathematica} yields that
 $$\sum_{k=0}^\infty\f{\bi{2k}k^2}{32^k}=\f{\Gamma(1/4)}{\sqrt{2\pi}\,\Gamma(3/4)}=\f{\Gamma(1/4)^2}{2\pi\sqrt{\pi}}.$$
 For any odd prime $p$, Z.-H. Sun \cite{SZH} proved that
 $$\sum_{k=0}^{p-1}\f{\bi{2k}k^2}{32^k}\eq\begin{cases} 2x-p/(2x)\pmod{p^2}&\t{if}\ p=x^2+y^2\ (x,y\in\Z)\ \t{with}
 \ 4\mid x-1,\\0\pmod{p^2}&\t{if}\ p\eq3\pmod4,\end{cases}$$
 which was previously conjectured by the author (cf. \cite[Conjecture 5.5]{S11}).
 \end{remark}

 \begin{conjecture} [{\rm 2022-12-30}] We have
 \begin{equation}\sum_{k=0}^\infty\f{\bi{2k}k\bi{3k}k}{(-216)^k}(3H_{3k}-H_k)
 =\l(\log \f89\r)\sum_{k=0}^\infty\f{\bi{2k}k\bi{3k}k}{(-216)^k}.
 \end{equation}
 \end{conjecture}
 \begin{remark}
 For any prime $p>3$,  we have
 $$\sum_{k=0}^{p-1}\f{\bi{2k}k\bi{3k}k}{(-216)^k}\eq\l(\f p3\r)\sum_{k=0}^{p-1}\f{\bi{2k}k\bi{3k}k}{24^k}\pmod{p^2}$$
 by \cite[Corollary 1.4]{FFA13}, and
 $$\sum_{k=0}^{p-1}\f{\bi{2k}k\bi{3k}k}{24^k}\eq\begin{cases}
 \bi{(2p-2)/3}{(p-1)/3}\pmod{p^2}&\t{if}\ p\eq1\pmod3,
 \\p/\bi{(2p+2)/3}{(p+1)/3}\pmod{p^2}&\t{if}\ p\eq2\pmod3,\end{cases}
 $$
 as conjectured by the author \cite[Conjecture 5.13]{S11c} and proved by C. Wang and Sun \cite[Theorem 1.2]{WS}.
 \end{remark}

 \begin{conjecture} [{\rm 2023-09-12}] {\rm (i)} We have
\begin{equation}\sum_{k=1}^\infty\f{48^k}{k(2k-1)\bi{2k}k\bi{4k}{2k}}
\l(6H_{4k-1}-9H_{2k-1}+2H_{k-1}+\f6{2k-1}\r)=\f{2\pi^3}{\sqrt3}.
\end{equation}

{\rm (ii)} Let $p>3$ be a prime. Then
\begin{equation}\sum_{k=1}^{(p-1)/2}\f{k\bi{2k}k\bi{4k}{2k}}{(2k+1)48^k}\l(6H_{4k}-9H_{2k}+2H_k-\f6{2k+1}\r)\eq0\pmod p\end{equation}
and
\begin{equation}\sum_{k=1}^{p-1}\f{k\bi{2k}k\bi{4k}{2k}}{(2k+1)48^k}\l(6H_{4k}-9H_{2k}+2H_k-\f6{2k+1}\r)
\eq\f{5p}{12}B_{p-2}\l(\f13\r)
\pmod {p^2}.\end{equation}
\end{conjecture}
\begin{remark} The author's conjectural identity
$$\sum_{k=1}^\infty\f{48^k}{k(2k-1)\bi{2k}k\bi{4k}{2k}}=\f{15}2K$$
(cf. \cite{S15}) was confirmed by Au \cite{Au22}.
\end{remark}

 \begin{conjecture} [{\rm 2022-12-30}] We have
 \begin{equation}\sum_{k=0}^\infty\f{\bi{4k}{2k}\bi{2k}k}{128^k}(2H_{4k}-H_{2k})=\f{(\log 2)\sqrt{\pi}}{2\Gamma(5/8)\Gamma(7/8)}.
 \end{equation}
 \end{conjecture}
 \begin{remark} {\tt Mathematica} yields the identity
 $$\sum_{k=0}^\infty\f{\bi{4k}{2k}\bi{2k}k}{128^k}=\f{\sqrt{\pi}}{\Gamma(5/8)\Gamma(7/8)}.$$
 By \cite[Corollary 1.3]{FFA13}, for any prime $p\eq5,7\pmod8$ we have
 $$\sum_{k=0}^{p-1}\f{\bi{4k}{2k}\bi{2k}k}{128^k}\eq0\pmod{p^2}.$$
 \end{remark}

 \begin{conjecture} [{\rm 2022-12-30}] We have
 \begin{equation}\sum_{k=0}^\infty\f{\bi{4k}{2k}\bi{2k}k}{72^k}(2H_{4k}-H_{2k})=(\log 3)\sum_{k=0}^\infty\f{\bi{4k}{2k}\bi{2k}k}{72^k}
 \end{equation}
 and
 \begin{equation}\sum_{k=0}^\infty\f{\bi{4k}{2k}\bi{2k}k}{576^k}(2H_{4k}-H_{2k})=\f12\l(\log \f 98\r)\sum_{k=0}^\infty\f{\bi{4k}{2k}\bi{2k}k}{576^k}.
 \end{equation}
 \end{conjecture}
 \begin{remark} Let $p>3$ be a prime. By \cite[Corollary 1.4]{FFA13},
 $$\sum_{k=0}^{p-1}\f{\bi{4k}{2k}\bi{2k}k}{576^k}\eq\l(\f{-2}p\r)\sum_{k=0}^{p-1}\f{\bi{4k}{2k}\bi{2k}k}{72^k}
 \pmod{p^2}.$$
 We also have
 \begin{align*}&\l(\f 6p\r)\sum_{k=0}^{p-1}\f{\bi{4k}{2k}\bi{2k}k}{72^k}
 \\\eq\ &\begin{cases}2x-\f p{2x}\pmod{p^2}&\t{if}\ p=x^2+y^2\ (x,y\in\Z\ \&\ 4\mid x-1),
 \\\f{2p}{3\bi{(p+1)/2}{(p+1)/4}}\pmod{p^2}&\t{if}\ p\eq3\pmod4,\end{cases}
 \end{align*}
 as conjectured by the author \cite[Conjecture 5.14(iii)]{S11c} and proved by Wang and Sun \cite[Theorem 5.2 and Remark 5.2]{WS}.
 \end{remark}

 \begin{conjecture} [{\rm 2022-12-30}] We have
 \begin{equation}\sum_{k=0}^\infty\f{\bi{4k}{2k}\bi{2k}k}{(-192)^k}(2H_{4k}-H_{2k})=\f12\l(\log \f34\r)\sum_{k=0}^\infty
 \f{\bi{4k}{2k}\bi{2k}k}{(-192)^k}.
 \end{equation}
 \end{conjecture}
 \begin{remark} Let $p>3$ be a prime. By \cite[Corollary 1.4]{FFA13}, we have
 $$\sum_{k=0}^{p-1}\f{\bi{4k}{2k}\bi{2k}k}{(-192)^k}\eq\l(\f{-2}p\r)\sum_{k=0}^{p-1}\f{\bi{4k}{2k}\bi{2k}k}{48^k}
 \pmod{p^2}.$$
 If $p=x^2+3y^2$ with $x,y\in\Z$ and $x\eq1\pmod3$,
 then
 $$\sum_{k=0}^{p-1}\f{\bi{4k}{2k}\bi{2k}k}{48^k}\eq2x-\f p{2x}\pmod{p^2},$$
 as conjectured by the author \cite[Conjecture 5.14]{S11c} and confirmed by G.-S. Mao and H. Pan \cite{MP}.
 If $p\eq2\pmod3$, then
 $$\sum_{k=0}^{p-1}\f{\bi{4k}{2k}\bi{2k}k}{48^k}\eq\f{3p}{2\bi{(p+1)/2}{(p+1)/6}}\pmod{p^2}$$
 as conjectured by the author \cite[Conjecture 5.14]{S11c} and confirmed by Wang and Sun \cite{WS}.
 \end{remark}

 \begin{conjecture} [{\rm 2022-12-30}] We have
 \begin{equation}\sum_{k=0}^\infty\f{\bi{4k}{2k}\bi{2k}k}{(-4032)^k}(2H_{4k}-H_{2k})=\f12\l(\log \f {63}{64}\r)\sum_{k=0}^\infty\f{\bi{4k}{2k}\bi{2k}k}{(-4032)^k}.
 \end{equation}
 \end{conjecture}
 \begin{remark} Let $p>3$ be a prime with $p\not=7$. By \cite[Corollary 1.4]{FFA13},
 $$\sum_{k=0}^{p-1}\f{\bi{4k}{2k}\bi{2k}k}{(-4032)^k}
 \eq\l(\f{-2}p\r)\sum_{k=0}^{p-1}\f{\bi{4k}{2k}\bi{2k}k}{63^k}\pmod{p^2}.$$
 The author \cite[Conjecture 5.14(ii)]{S11c} conjectured that
 \begin{align*}&\sum_{k=0}^{p-1}\f{\bi{4k}{2k}\bi{2k}k}{63^k}
 \\\eq\ &\begin{cases}(\f p3)(2x-\f p{2x})\pmod{p^2}&\t{if}\ p=x^2+7y^2\ \t{with}\ x,y\in\Z\ \t{and}\ (\f x7)=1,\\0\pmod p&\t{if}\ (\f p7)=-1.\end{cases}
 \end{align*}
 \end{remark}

 \begin{conjecture} [{\rm 2022-12-30}] We have
 \begin{equation}\sum_{k=0}^\infty\f{\bi{6k}{3k}\bi{3k}k}{864^k}
 (6H_{6k}-3H_{3k}-2H_{2k}+H_k)=\f{(\log 2)\sqrt{\pi}}{\Gamma(7/12)\Gamma(11/12)}.
 \end{equation}
 \end{conjecture}
 \begin{remark} {\tt Mathematica} yields the identity
 $$\sum_{k=0}^\infty\f{\bi{6k}{3k}\bi{3k}k}{864^k}=\f{\sqrt{\pi}}{\Gamma(7/12)\Gamma(11/12)}.$$
 By \cite[Corollary 1.3]{FFA13}, for any prime $p>3$ with $p\eq3\pmod4$ we have
 $$\sum_{k=0}^{p-1}\f{\bi{6k}{3k}\bi{3k}k}{864^k}\eq0\pmod{p^2}.$$
 \end{remark}

 \section{Series and congruences with summands \\containing $3$ or $4$ binomial coefficients}

 \setcounter{theorem}{0}
 \setcounter{proposition}{0}

 In 1993 D. Zeilberger \cite{Z} used the WZ method to establish the identity
 $$\sum_{k=1}^\infty\f{21k-8}{k^3\bi{2k}k^3}=\f{\pi^2}6.$$
 The author \cite{S15i} proved that
 $$\sum_{k=1}^{(p-1)/2}\f{21k-8}{k^3\bi{2k}k^3}\eq-\l(\f{-1}p\r)4E_{p-3}\pmod p$$
 for any prime $p>3$.

 \begin{conjecture} [{\rm 2022-10-11}] \label{21Z} {\rm (i)} We have the identity
 \begin{equation}\label{2880}
 \sum_{k=1}^\infty\f{21k-8}{k^3\bi{2k}k^3}\l(H_{2k-1}^{(2)}-\f{25}8H_{k-1}^{(2)}\r)=\f{47\pi^4}{2880}.
 \end{equation}

 {\rm (ii)} For any prime $p>3$, we have
 \begin{equation}
 \sum_{k=1}^{(p-1)/2}(21k+8)\bi{2k}k^3\l(H_{2k}^{(2)}-\f{25}8H_{k}^{(2)}\r)\eq32p\l(\f{-1}p\r)E_{p-3}
 \pmod{p^2}
 \end{equation}
 and
 \begin{equation}
 \sum_{k=0}^{p-1}(21k+8)\bi{2k}k^3\l(H_{2k}^{(2)}-\f{25}8H_{k}^{(2)}\r)\eq-48H_{p-1}
 +\f{246}5p^4B_{p-5} \pmod{p^5}.
 \end{equation}

 {\rm (iii)} For each prime $p>5$, we have
\begin{equation}\sum_{k=1}^{p-1}\bi{2k}k^3\l(2(21k+8)(H_{2k}-H_k)+7\r)\eq-112 pH_{p-1}\pmod{p^5}.
\end{equation}
 \end{conjecture}
 \begin{remark} After seeing an earlier arXiv version of this paper, Au \cite[Corollary 2.3]{Au22}
 confirmed the author's conjectural identity \eqref{2880},
 and proved an identity (after the proof of \cite[Theorem 2.2]{Au22}) which has the equivalent form
 $$\sum_{k=1}^\infty\f{(21k-8)(H_{2k-1}-H_{k-1})-7/2}{k^3\bi{2k}k^3}=\zeta(3).$$
 \end{remark}

 \begin{conjecture} [{\rm 2022-10-11}] \label{21cube} {\rm (i)} We have
 \begin{equation} \label{21H3}\sum_{k=1}^\infty\f{21k-8}{k^3\bi{2k}k^3}\l(H_{2k-1}^{(3)}+\f{43}8H_{k-1}^{(3)}\r)
 =\f{711}{28}\zeta(5)-\f{29}{14}\pi^2\zeta(3).
 \end{equation}

 {\rm (ii)} For any prime $p>7$, we have
 \begin{equation} \sum_{k=0}^{(p-1)/2}(21k+8)\bi{2k}k^3\l(H_{2k}^{(3)}+\f{43}8H_k^{(3)}\r)
 \eq32\l(\f{-1}p\r)E_{p-3}\pmod p
 \end{equation}
 and
 \begin{equation} \sum_{k=0}^{p-1}(21k+8)\bi{2k}k^3\l(H_{2k}^{(3)}+\f{43}8H_k^{(3)}\r)
 \eq-\f{120}7pB_{p-3}\pmod {p^2}.
 \end{equation}
 \end{conjecture}
 \begin{remark} The identity \eqref{21H3} looks quite challenging.
 \end{remark}

 \begin{conjecture} [{\rm 2022-10-13}] {\rm (i)} We have
 \begin{equation}\label{beta}\sum_{k=1}^{\infty}\f{(3k-1)(-8)^k}{k^3\bi{2k}k^3}\l(H_{2k-1}^{(2)}-\f54H_{k-1}^{(2)}\r)
 =-2\beta(4).
 \end{equation}

 {\rm (ii)} For any prime $p>3$, we have
 \begin{equation}\sum_{k=0}^{(p-1)/2}(3k+1)\f{\bi{2k}k^3}{(-8)^k}\l(H_{2k}^{(2)}-\f54H_k^{(2)}\r)
 \eq\l(\f 2p\r)\f p4E_{p-3}\l(\f14\r)\pmod {p^2}
 \end{equation}
 and
 \begin{equation}\sum_{k=0}^{p-1}(3k+1)\f{\bi{2k}k^3}{(-8)^k}\l(H_{2k}^{(2)}-\f54H_k^{(2)}\r)
 \eq pE_{p-3}\pmod {p^2}.
 \end{equation}

 {\rm (iii)} Let $p$ be an odd prime. Then
 \begin{equation}\sum_{k=0}^{(p-1)/2}\f{\bi{2k}k^3}{(-8)^k}
 \l(2(3k+1)(H_{2k}-H_k)+1\r)\eq\l(\f{-1}p\r)2^{p-1}\pmod{p^2}
 \end{equation}
 and
 \begin{equation}\sum_{k=0}^{p-1}(3k+1)\f{\bi{2k}k^3}{(-8)^k}
 \l(H_{2k}^{(3)}+\f 78H_k^{(3)}\r)\eq0\pmod p.
 \end{equation}
 \end{conjecture}
\begin{remark} In 2008, J. Guillera \cite{G} used the WZ method to find the identity
$$\sum_{k=1}^\infty\f{(3k-1)(-8)^k}{k^3\bi{2k}k^3}=-2G.$$
The identity \eqref{beta} provides a fast converging series for computing the constant $\beta(4)$.
We are unable to find the exact values of the series
$$\sum_{k=1}^\infty\f{(-8)^k}{k^3\bi{2k}k^3}(2(3k-1)(H_{2k-1}-H_{k-1})-1)$$
and
$$\sum_{k=1}^{\infty}\f{(3k-1)(-8)^k}{k^3\bi{2k}k^3}
 \l(H_{2k-1}^{(3)}+\f 78H_{k-1}^{(3)}\r).$$
\end{remark}

 \begin{conjecture}[{\rm 2022-10-11}] \label{16} {\rm (i)} We have the identity
 \begin{equation}\label{Pi^4}
 \sum_{k=1}^\infty\f{(3k-1)16^k}{k^3\bi{2k}k^3}\l(H_{2k-1}^{(2)}-\f{5}4H_{k-1}^{(2)}\r)=\f{\pi^4}{24}.
 \end{equation}

 {\rm (ii)} Let $p>3$ be a prime. Then
 \begin{equation}
 \sum_{k=0}^{(p-1)/2}(3k+1)\f{\bi{2k}k^3}{16^k}\l(H_{2k}^{(2)}-\f{5}4H_{k}^{(2)}\r)
 \eq2p\l(\f{-1}p\r)E_{p-3}
 \pmod{p^2}
 \end{equation}
 and
 \begin{equation}\label{GL-c}
 \sum_{k=0}^{p-1}(3k+1)\f{\bi{2k}k^3}{16^k}\l(H_{2k}^{(2)}-\f{5}4H_{k}^{(2)}\r)
 \eq\f 76p^2B_{p-3}
 \pmod{p^3}.
 \end{equation}

 {\rm (iii)} Let $p$ be an odd prime. Then
 \begin{equation}
 \sum_{k=1}^{p-1}\f{\bi{2k}k^3}{16^k}\l(2(3k+1)(H_{2k}-H_k)+1\r)
 \eq\f43p\,q_p(2)-\f23p^2q_p(2)^2 \pmod{p^3}.
 \end{equation}
 If $p\eq2\pmod3$, then
 \begin{equation}
 \sum_{k=0}^{(p-1)/2}\f{\bi{2k}k^3}{16^k}\l(H_{2k}^{(2)}-\f{5}4H_{k}^{(2)}\r)
 \eq0\pmod{p}.
 \end{equation}
 \end{conjecture}
 \begin{remark} In 2008 Guillera \cite[Identity 1]{G} used the WZ method to establish the identity
 $$\sum_{k=1}^\infty\f{(3k-1)16^k}{k^3\bi{2k}k^3}=\f{\pi^2}2.$$
 Two $q$-analogues of this identity were given by Q.-H. Hou, C. Krattenthaler and Z.-W. Sun \cite{HKS}.
 Guo and Lian \cite{GL} proved that the two sides of \eqref{GL-c} are congruent modulo $p^2$ for any prime $p>3$.
 After seeing an earlier arXiv version of this paper, Au \cite{Au22}
 confirmed the author's conjectural identity \eqref{Pi^4},
 and proved an identity (after the proof of \cite[Corollary 2.3]{Au22}) which has the equivalent form
 $$\sum_{k=1}^\infty\f{16^k}{k^3\bi{2k}k^3}\l((3k-1)(H_{2k-1}-H_{k-1})-\f12\r)
 =\f{\pi^2}3\log2+\f 76\zeta(3).$$
 \end{remark}

 \begin{conjecture} [{\rm 2022-10-11}] \label{16cube} {\rm (i)} We have
 \begin{equation}\sum_{k=1}^\infty\f{(3k-1)16^k}{k^3\bi{2k}k^3}\l(H_{2k-1}^{(3)}+\f78H_{k-1}^{(3)}\r)
 =\f{\pi^2}2\zeta(3).
 \end{equation}

 {\rm (ii)} For any odd prime $p$, we have
 \begin{equation}\sum_{k=1}^{(p-1)/2}(3k+1)\f{\bi{2k}k^3}{16^k}\l(H_{2k}^{(3)}+\f78H_{k}^{(3)}\r)
 \eq 2\l(\f{-1}p\r)E_{p-3}\pmod p
 \end{equation}
 and
 \begin{equation}\sum_{k=0}^{p-1}(3k+1)\f{\bi{2k}k^3}{16^k}\l(H_{2k}^{(3)}+\f78H_{k}^{(3)}\r)
 \eq \f32pB_{p-3}\pmod {p^2}.
 \end{equation}
 \end{conjecture}
 \begin{remark} Conjecture \ref{16cube} looks more challenging than Conjecture \ref{16}.
 \end{remark}

 \begin{conjecture} [{\rm 2022-10-16}] {\rm (i)} We have
 \begin{equation}\label{16beta}\sum_{k=1}^\infty\f{(4k-1)(-64)^k}{k^3\bi{2k}k^3}\l(H_{2k-1}^{(2)}-\f12H_{k-1}^{(2)}\r)
 =-16\beta(4).
 \end{equation}

 {\rm (ii)} For any prime $p>3$, we have
 \begin{align}
 \label{-64}\sum_{k=0}^{(p-1)/2}(4k+1)\f{\bi{2k}k^3}{(-64)^k}\l(H_{2k}^{(2)}-\f1{2}H_k^{(2)}\r)
 &\eq pE_{p-3}\pmod{p^2},
 \\\label{-64+}\sum_{k=0}^{p-1}(4k+1)\f{\bi{2k}k^3}{(-64)^k}\l(H_{2k}^{(2)}-\f1{2}H_k^{(2)}\r)
 &\eq pE_{p-3}\pmod{p^2}.
 \end{align}
 \end{conjecture}
 \begin{remark} Guillera \cite{G} proved that
 $$\sum_{k=1}^\infty\f{(4k-1)(-64)^k}{k^3\bi{2k}k^3}=-16G.$$
 The congruence \eqref{-64} is motivated by \eqref{64}. After seeing an earlier arXiv version of this paper,  Au \cite[Corollary 2.11]{Au22}
 confirmed the author's conjectural identity \eqref{16beta}.
 It seems that for any $m,n\in\Z^+$ and odd prime $p$, we have
\begin{equation}\sum_{k=0}^{(p-1)/2}(4k+1)\f{\bi{2k}k^n}{(-4)^{kn}}\l(H_{2k}^{(2m)}-\f{H_k^{(2m)}}{2^{2m-1}}\r)\eq0\pmod p.\end{equation}
 \end{remark}

 \begin{conjecture} [{\rm 2022-12-05}] {\rm (i)} We have
 \begin{equation}
 \sum_{k=1}^\infty\f{8^{k}((10k-3)(H_{2k-1}-H_{k-1})-1)}{k^3\bi{2k}k^2\bi{3k}k}=\f72\zeta(3)
 \end{equation}
 and
  \begin{equation}
 \sum_{k=1}^\infty\f{8^{k}((10k-3)(H_{3k-1}-H_{k-1})-8/3)}{k^3\bi{2k}k^2\bi{3k}k}=\f{2\pi^2\log2+7\zeta(3)}4.
 \end{equation}

 {\rm (ii)} For any odd prime $p$, we have
 \begin{equation}\sum_{k=1}^{p-1}\f{\bi{2k}k^2\bi{3k}k}{8^k}\l((10k+3)(H_{2k}-H_k)+1\r)\eq\f{63}8p^3B_{p-3}
 \pmod{p^4}
 \end{equation}
 and
 \begin{equation}\sum_{k=1}^{p-1}\f{\bi{2k}k^2\bi{3k}k}{8^k}
 \l(3(10k+3)(H_{3k}-H_k)+8\r)
 \eq9p\,q_p(2)-\f 92p^2q_p(2)^2\pmod{p^3}.
 \end{equation}
 \end{conjecture}
 \begin{remark} As conjectured by the author \cite{S11c} and confirmed by Guillera and M. Rogers \cite{GR}, we have
 $$\sum_{k=1}^\infty\f{(10k-3)8^{k}}{k^3\bi{2k}k^2\bi{3k}k}=\f{\pi^2}2.$$
 \end{remark}

 \begin{conjecture} [{\rm 2022-12-05}]\label{27^k} {\rm (i)} We have
 \begin{equation}
 \sum_{k=1}^\infty\f{(-27)^{k}((15k-4)(3H_{3k-1}-H_{k-1})-9)}{k^3\bi{2k}k^2\bi{3k}k}=
 -\f{4\pi^3}{\sqrt3}.
 \end{equation}

 {\rm (ii)} For any prime $p>3$, we have
 \begin{equation}\sum_{k=0}^{p-1}\f{\bi{2k}k^2\bi{3k}k}{(-27)^k}((15k+4)(3H_{3k}-H_k)+9)
 \eq9\l(\f p3\r)+6p^2B_{p-2}\l(\f 13\r)\pmod{p^3}.
 \end{equation}
 \end{conjecture}
 \begin{remark} As conjectured by the author \cite{S11c} and confirmed by Kh. Hessami Pilehrood and T. Hessami Pilehrood \cite{HP}, we have
 $$\sum_{k=1}^\infty\f{(15k-4)(-27)^{k-1}}{k^3\bi{2k}k^2\bi{3k}k}=K.$$
 \end{remark}

 \begin{conjecture}[{\rm 2023-08-26}]
{\rm (i)} We have
\begin{equation} \sum_{k=1}^\infty\f{(-27)^k}{k^3\bi{2k}k^2\bi{3k}k}\l((15k-4)(H_{2k-1}-H_{k-1})+\f1{2k-1}\r)
=-\f{4\pi^3}{3\sqrt3}.
\end{equation}

{\rm (ii)} For any prime $p>3$, we have
\begin{equation}\begin{aligned}&\ \sum_{k=0}^{p-1}\f{\bi{2k}k^2\bi{3k}k}{(-27)^k}\l((15k+4)(H_{2k}-H_k)+\f1{2k+1}\r)
\\&\quad\eq\l(\f p3\r)+\f43p^2B_{p-2}\l(\f13\r)\pmod{p^3}.
\end{aligned}
\end{equation}
\end{conjecture}
\begin{remark} This conjecture is similar to Conjecture \ref{27^k}.
\end{remark}

 \begin{conjecture} [{\rm 2022-12-05}] {\rm (i)} We have
 \begin{equation}
 \sum_{k=1}^\infty\f{64^{k-1}((11k-3)(2H_{2k-1}+H_{k-1})-4)}{k^3\bi{2k}k^2\bi{3k}k}=\f72\zeta(3)
 \end{equation}
 and
 \begin{equation}
 \sum_{k=1}^\infty\f{64^{k-1}((11k-3)(3H_{3k-1}-6H_{k-1})-7)}{k^3\bi{2k}k^2\bi{3k}k}=\f{6\pi^2\log2-21\zeta(3)}8.
 \end{equation}

 {\rm (ii)} For any odd prime $p$, we have
 \begin{equation}\sum_{k=1}^{p-1}\f{\bi{2k}k^2\bi{3k}k}{64^k}((11k+3)(2H_{2k}+H_k)+4)
 \eq21p^3B_{p-3}\pmod{p^4}
 \end{equation}
 and
 \begin{equation}\sum_{k=1}^{p-1}\f{\bi{2k}k^2\bi{3k}k}{64^k}((11k+3)(3H_{3k}-6H_k)+7)
 \eq18p\,q_p(2)+9p^2q_p(2)^2\pmod{p^3}.
 \end{equation}
 \end{conjecture}
 \begin{remark} As conjectured by the author \cite{S11c} and confirmed by Guillera \cite{Gu}, we have
 $$\sum_{k=1}^\infty\f{(11k-3)64^k}{k^3\bi{2k}k^2\bi{3k}k}=8\pi^2.$$
 \end{remark}

 \begin{conjecture} [{\rm 2022-12-09}] {\rm (i)} We have
 \begin{equation}
 \sum_{k=1}^\infty\f{81^{k}((35k-8)(H_{4k-1}-H_{k-1})-35/4)}{k^3\bi{2k}k^2\bi{4k}{2k}}=12\pi^2\log3+39\zeta(3).
 \end{equation}

 {\rm (ii)} For any prime $p>3$, we have
 \begin{equation}\begin{aligned}&\sum_{k=1}^{p-1}\f{\bi{2k}k^2\bi{4k}{2k}}{81^k}
 \l(4(35k+8)(H_{4k}-H_k)+35\r)
 \\\eq\ &32(3^{p-1}-1)-16(3^{p-1}-1)^2\pmod{p^3}.
 \end{aligned}
 \end{equation}
 \end{conjecture}
 \begin{remark} As conjectured by the author \cite{S11c} and confirmed in \cite{GR}, we have
 $$\sum_{k=1}^\infty\f{(35k-8)81^k}{k^3\bi{2k}k^2\bi{4k}{2k}}=12\pi^2.$$
 \end{remark}

\begin{conjecture} [{\rm 2023-01-14}] For any prime $p>3$, we have
 \begin{equation}\begin{aligned}&\sum_{k=1}^{p-1}\f{\bi{2k}k^2\bi{4k}{2k}}{(-144)^k}
 \l(4(5k+1)(H_{4k}-H_k)+5\r)
 \\\eq\ &\l(\f p3\r)\l(5+2p(2q_p(2)+q_p(3))\r)\pmod{p^2}.
 \end{aligned}
 \end{equation}
 \end{conjecture}
 \begin{remark} As conjectured by the author \cite{S11c} and confirmed in \cite{GR}, we have
 $$\sum_{k=1}^\infty\f{(5k-1)(-144)^k}{k^3\bi{2k}k^2\bi{4k}{2k}}=-\f{45}2K.$$
 We are unable to find the exact value of the series
 $$\sum_{k=1}^\infty\f{(-144)^k}{k^3\bi{2k}k^2\bi{4k}{2k}}(4(5k-1)(H_{4k-1}-H_{k-1})-5).$$
 \end{remark}

 The classical rational Ramanujan-type series for $1/\pi$ have the following four forms:
\begin{align}\sum_{k=0}^\infty(ak+b)\frac{\binom{2k}k^3}{m^k}&=\frac{c\sqrt d}{\pi},\label{1}
\\\sum_{k=0}^\infty(ak+b)\frac{\binom{2k}k^2\binom{3k}k}{m^k}&=\frac{c\sqrt d}{\pi},\label{2}
\\\sum_{k=0}^\infty(ak+b)\frac{\binom{2k}k^2\binom{4k}{2k}}{m^k}&=\frac{c\sqrt d}{\pi},\label{3}
\\\sum_{k=0}^\infty(ak+b)\frac{\binom{2k}k\binom{3k}k\binom{6k}{3k}}{m^k}&=\frac{c\sqrt d}{\pi},\label{4}
\end{align}
where $a,b,m\in\mathbb Z$, $am\not=0$, $c\in\mathbb Q\sm\{0\}$, and $d$ is a positive squarefree integer. It is known that there are totally $36$ such series, see, e.g., S. Cooper \cite[Chapter 14]{Co17}.

For a positive integer $m$, can we find similar series for $(\log m)/\pi$? Motivated by Ramanujan-type series of the forms \eqref{1}-\eqref{4}, the author formulated the following general conjecture.

\begin{conjecture} [{\rm General Conjecture, 2022-12-08}] \label{abm} {\rm (i)} If we have an identity \eqref{1} with $a,b,m\in\mathbb Z$, $am\not=0$, $c\in\mathbb Q\sm\{0\}$, and $d\in\Z^+$ squarefree, then
\begin{equation}\label{I}\sum_{k=0}^\infty\frac{\binom{2k}k^3}{m^k}(6(ak+b)(H_{2k}-H_k)+a)
=\f{c\sqrt d}{\pi}\log|m|,
\end{equation}
and
\begin{equation}\sum_{k=0}^{p-1}\f{\bi{2k}k^3}{m^k}\l(6(ak+b)(H_{2k}-H_k)+a\r)\eq \l(\f{-d}p\r)(a+b(m^{p-1}-1))\pmod{p^2}
\end{equation}
for any odd prime $p\nmid dm$.

{\rm (ii)} If we have an identity \eqref{2} with $a,b,m\in\mathbb Z$, $am\not=0$, $c\in\mathbb Q\sm\{0\}$, and $d\in\Z^+$ squarefree, then
\begin{equation}\label{II}\sum_{k=0}^\infty\frac{\binom{2k}k^2\binom{3k}k}{m^k}((ak+b)(3H_{3k}+2H_{2k}-5H_k)+a)
=\f{c\sqrt d}{\pi}\log|m|,
\end{equation}
and
\begin{equation}\begin{aligned}&\sum_{k=0}^{p-1}\f{\bi{2k}k^2\bi{3k}k}{m^k}
\l((ak+b)(3H_{3k}+2H_{2k}-5H_k)+a\r)
\\&\qquad\eq\l(\f{-d}p\r)(a+b(m^{p-1}-1))\pmod{p^2}
\end{aligned}
\end{equation}
for any odd prime $p\nmid dm$.

{\rm (iii)} If we have an identity \eqref{3} with $a,b,m\in\mathbb Z$, $am\not=0$, $c\in\mathbb Q\sm\{0\}$,
and $d\in\Z^+$ squarefree, then
\begin{equation}\label{III}\sum_{k=0}^\infty\frac{\binom{2k}k^2\binom{4k}{2k}}{m^k}(4(ak+b)(H_{4k}-H_k)+a)
=\f{c\sqrt d}{\pi}\log|m|,
\end{equation}
and
\begin{equation}\label{III'}
\begin{aligned}&\sum_{k=0}^{p-1}\frac{\binom{2k}k^2\binom{4k}{2k}}{m^k}(4(ak+b)(H_{4k}-H_k)+a)
\\&\qquad\eq \l(\f{-d}p\r)\l(a+b(m^{p-1}-1)\r)\pmod{p^2}
\end{aligned}
\end{equation}
for any odd prime $p\nmid dm$.

{\rm (iv)} If we have an identity \eqref{4} with $a,b,m\in\mathbb Z$, $am\not=0$, $c\in\mathbb Q\sm\{0\}$, and $d\in\Z^+$ squarefree, then
\begin{equation}\label{IV}\sum_{k=0}^\infty\frac{\binom{2k}k\binom{3k}{k}\binom{6k}{3k}}{m^k}(3(ak+b)
(2H_{6k}-H_{3k}-H_k)+a)=\f{c\sqrt d}{\pi}\log|m|,
\end{equation}
and
\begin{equation}\label{IV'}\begin{aligned}&\sum_{k=0}^{p-1}\frac{\binom{2k}k\binom{3k}{k}\binom{6k}{3k}}{m^k}(3(ak+b)
(2H_{6k}-H_{3k}-H_k)+a)
\\&\qquad\eq\l(\f{-d}p\r)\l(a+b(m^{p-1}-1)\r)\pmod{p^2}
\end{aligned}
\end{equation}
for any odd prime $p\nmid dm$.
\end{conjecture}
\begin{remark} Ramanujan \cite{R} found the irrational series
$$\sum_{k=0}^\infty\left(k+\frac{31}{270+48\sqrt5}\right)\frac{\binom{2k}k^3}{(2^{20}/(\sqrt5-1)^8)^k}=\frac{16}{(15+21\sqrt5)\pi}.$$
In the spirit of part (i) of our general conjecture, we guess that
\begin{align*}&\sum_{k=0}^\infty\frac{\binom{2k}k^3}{(2^{20}/(\sqrt5-1)^8)^k}\left(6\left(k+\frac{31}{270+48\sqrt5}\right)(H_{2k}-H_k)+1\right)
\\&\qquad\ = \frac{16}{(15+21\sqrt5)\pi}\times\log\frac{2^{20}}{(\sqrt5-1)^8},
\end{align*}
which can be easily checked via {\tt Mathematica}. We believe that similar things happen
for all irrational Ramanujan-type series.
\end{remark}

 {\it Example}\ 3.1.
 In view of the Ramanujan series (cf. \cite[Chapter 14]{Co17} and \cite{R})
 $$\sum_{k=0}^\infty(5k+1)\f{\bi{2k}k^2\bi{3k}k}{(-192)^k}=\f{4\sqrt3}{\pi}$$
 and the known identity
 $$\sum_{k=0}^\infty\f{\bi{2k}k^2\bi{3k}k}{(-192)^k}=\f{4\sqrt{\pi}}{3\sqrt3\,\Gamma(5/6)^3}$$
 (cf. (14.29) of \cite[p.\,624]{Co17}), by part (ii) of Conjecture \ref{abm} we should have
 \begin{equation*}\label{192Hk}
 \sum_{k=0}^\infty\f{\bi{2k}k^2\bi{3k}k}{(-192)^k}(5k+1)(3H_{3k}+2H_{2k}-5H_k)
 =\f{4\log192}{\sqrt3\,\pi}-\f{20\sqrt{\pi}}{3\sqrt3\,\Gamma(5/6)^3}.
 \end{equation*}

 In 2013, Guillera \cite[(32)]{G13} proved the identity
 $$\sum_{k=0}^\infty\f{\bi{2k}k^3}{(-64)^k}\l((4k+1)H_k-\f23\r)=-\f{4\log2}{\pi}$$
 motivated by the Bauer series \eqref{Bauer}.

 \begin{conjecture} [{\rm 2022-10-16}] We have
 \begin{equation}\sum_{k=0}^\infty(4k+1)\f{\bi{2k}k^3}{(-64)^k}H_{2k}^{(3)}=\f{15\zeta(3)}{4\pi}-2G.
 \end{equation}
 \end{conjecture}
 \begin{remark}
 For any $m,n\in\Z^+$ and odd prime $p$ not dividing $2^{2m-1}-1$, we have
 $H_{p-1}^{(2m-1)}\eq0\pmod p$ since
 $\sum_{j=1}^{p-1}1/{(2j)^{2m-1}}\eq\sum_{k=1}^{p-1}1/{k^{2m-1}}\pmod p,$
 thus
 \begin{align*}&\sum_{k=0}^{(p-1)/2}(4k+1)\bi{(p-1)/2}k^nH_{2k}^{(2m-1)}
 \\=\ &\sum_{k=0}^{(p-1)/2}\l(4\l(\f{p-1}2-k\r)+1\r)\bi{(p-1)/2}k^nH_{p-1-2k}^{(2m-1)}
 \\\eq\ &-\sum_{k=0}^{(p-1)/2}\l(4k+1\r)\bi{(p-1)/2}k^nH_{2k}^{(2m-1)}\pmod p
 \end{align*}
 and hence
 \begin{equation}\sum_{k=0}^{(p-1)/2}(4k+1)\f{\bi{2k}k^n}{(-4)^{kn}}H_{2k}^{(2m-1)}\eq0\pmod p.
 \end{equation}
 (Note that $\bi{-1/2}k=\bi{2k}k/(-4)^k$ for any $k\in\N$.)
 \end{remark}

\begin{conjecture} [{\rm 2022-12-04}] {\rm (i)} We have
 \begin{equation}\label{Hk}
 \sum_{k=0}^\infty(6k+1)\f{\bi{2k}k^3}{256^k}H_k=\f{4}3\cdot\f{\root3\of2\sqrt{\pi}}{\Gamma(5/6)^3}-\f{8\log2}{\pi}
 \end{equation}
 and
 \begin{equation}\label{H2k}
 \sum_{k=0}^\infty(6k+1)\f{\bi{2k}k^3}{256^k}H_{2k}=\f{2}3\cdot\f{\root3\of2\sqrt{\pi}}{\Gamma(5/6)^3}-\f{8\log2}{3\pi}.
 \end{equation}

 {\rm (ii)} Let $p$ be an odd prime. Then
 \begin{equation}\sum_{k=0}^{(p-1)/2}\f{\bi{2k}k^3}{256^k}((6k+1)(3H_{2k}-H_k)-1)\eq(-1)^{(p+1)/2}\pmod{p^4}.
 \end{equation}
 If $p>3$, then
 \begin{equation}\begin{aligned}&\sum_{k=0}^{(p-1)/2}\f{\bi{2k}k^3}{256^k}((6k+1)(H_{2k}-H_k)+1)
 \\\eq\ &\l(\f{-1}p\r)\l(1+\f 43p\,q_p(2)-\f23p^2q_p(2)^2\r)\pmod{p^3}.
 \end{aligned}
 \end{equation}
 \end{conjecture}
 \begin{remark} It is known (cf. (14.27) of \cite[p.\,623]{Co17}) that
 $$\sum_{k=0}^\infty\f{\bi{2k}k^3}{256^k}=\f 23\cdot\f{\root3\of2\sqrt{\pi}}{\Gamma(5/6)^3}.$$
 In view of this,  \eqref{Hk} and \eqref{H2k} together implies the identities
 $$\sum_{k=0}^\infty\f{\bi{2k}k^3}{256^k}\l((6k+1)(3H_{2k}-H_k)-1\r)=0$$
 and $$\sum_{k=0}^\infty\f{\bi{2k}k^3}{256^k}\l((6k+1)(H_{2k}-H_k)+1\r)=\f{16\log2}{3\pi}.$$
 The last identity is also implied by Conjecture \ref{abm}(i).
 \end{remark}

 \begin{conjecture} [{\rm 2022-10-12}] Let $p>3$ be a prime.

{\rm (i)} We have
\begin{equation}\label{256-2}\sum_{k=0}^{(p-1)/2}(6k+1)\f{\bi{2k}k^3}{256^k}
 \l(H_{2k}^{(2)}-\f5{16}H_k^{(2)}\r)\eq\f 7{24}\l(\f{-1}p\r)p^2B_{p-3}\pmod{p^3},
 \end{equation}
 and
 \begin{equation}\sum_{k=0}^{p-1}(6k+1)\f{\bi{2k}k^3}{256^k}
 \l(H_{2k}^{(2)}-\f5{16}H_k^{(2)}\r)\eq-pE_{p-3}\pmod{p^2}.
 \end{equation}

{\rm (ii)} If $p\eq2\pmod3$, then
 \begin{equation}\sum_{k=0}^{(p-1)/2}\f{\bi{2k}k^3}{256^k}\l(H_{2k}^{(2)}-\f5{16}H_k^{(2)}\r)
 \eq0\pmod p.
 \end{equation}
 \end{conjecture}
 \begin{remark} For any prime $p>3$, Guo and Lian \cite{GL} proved that the two sides of
 \eqref{256-2} are congruent modulo $p^2$.
 \end{remark}

 \begin{conjecture} [{\rm 2022-10-11}] {\rm (i)} We have the identity
 \begin{equation}\sum_{k=0}^\infty (6k+1)\f{\bi{2k}k^3}{256^k}
 \l(H_{2k}^{(3)}-\f7{64}H_k^{(3)}\r)=\f{25\zeta(3)}{8\pi}-G.
 \end{equation}

 {\rm (ii)} Let $p$ be an odd prime. Then
 \begin{equation}\sum_{k=0}^{(p-1)/2}(6k+1)\f{\bi{2k}k^3}{256^k}
 \l(H_{2k}^{(3)}-\f7{64}H_k^{(3)}\r)\eq-\f 12E_{p-3}\pmod{p}
 \end{equation}
 and
 \begin{equation}\sum_{k=0}^{p-1}(6k+1)\f{\bi{2k}k^3}{256^k}
 \l(H_{2k}^{(3)}-\f7{64}H_k^{(3)}\r)\eq-\f 32E_{p-3}\pmod{p}.
 \end{equation}
 \end{conjecture}
 \begin{remark} For any $k\in\Z^+$, it is easy to see that
 $$H_{2k}^{(3)}-\f7{64}H_k^{(3)}=\sum_{j=1}^k\l(\f1{(2j-1)^3}+\f1{(4j)^3}\r).$$
 \end{remark}

 \begin{conjecture} [{\rm 2022-10-12}] Let $p>3$ be a prime. Then
 \begin{equation}\begin{aligned}\label{-512}\sum_{k=0}^{(p-1)/2}(6k+1)\f{\bi{2k}k^3}{(-512)^k}\l(H_{2k}^{(2)}-\f5{16}H_k^{(2)}\r)
 &\eq\l(\f 2p\r)\f p4E_{p-3}\pmod{p^2}
 \end{aligned}\end{equation}
 and
  \begin{equation}\begin{aligned}\label{-512+}\sum_{k=0}^{p-1}(6k+1)\f{\bi{2k}k^3}{(-512)^k}\l(H_{2k}^{(2)}-\f5{16}H_k^{(2)}\r)
 &\eq\f p{16}E_{p-3}\l(\f14\r)\pmod{p^2}.
 \end{aligned}\end{equation}
 \end{conjecture}
 \begin{remark} Note that \eqref{-512} is stronger than \eqref{L-cong}.
 \end{remark}

 \begin{conjecture} [{\rm 2022-10-16}] {\rm (i)} We have the identity
 \begin{equation}
 \sum_{k=0}^\infty(6k+1)\f{\bi{2k}k^3}{(-512)^k}\l(H_{2k}^{(3)}-\f{7}{64}H_{k}^{(3)}\r)
 =\f{57}{16}\cdot\f{\zeta(3)}{\sqrt2\,\pi}-L,
 \end{equation}
 where
 $$L:=L\l(2,\l(\f{-8}{\cdot}\r)\r)=\sum_{n=1}^\infty\f{(\f{-8}n)}{n^2}
 =\sum_{k=0}^\infty\f{(-1)^{k(k-1)/2}}{(2k+1)^2}$$
 with $(\f{-8}{\cdot})$ the Kronecker symbol.

 {\rm (ii)} Let $p$ be an odd prime. Then
 \begin{equation}\label{512-3-c}\sum_{k=0}^{p-1}(6k+1)\f{\bi{2k}k^3}{(-512)^k}\l(H_{2k}^{(3)}-\f 7{64}H_k^{(3)}\r)
 \eq0\pmod p.
 \end{equation}
 \end{conjecture}
 \begin{remark} For a general odd prime $p$, we are unable to find a closed form for the left-hand side of
the congruence \eqref{512-3-c} modulo $p^2$.
 \end{remark}

 \begin{conjecture}[{\rm 2022-12-09}] \label{2kk} For any prime $p>3$, we have
\begin{equation}\begin{aligned}&\sum_{k=0}^{(p-1)/2}\f{\bi{2k}k^3}{4096^k}\l((42k+5)(H_{2k}-H_k)+7\r)
\\\eq&\ \l(\f{-1}p\r)\l(7+10p\,q_p(2)-5p^2q_p(2)^2\r)\pmod{p^3}.
\end{aligned}
\end{equation}
\end{conjecture}
\begin{remark}\label{4096R} In view of the Ramanujan series (cf. \cite{R})
$$\sum_{k=0}^\infty(42k+5)\f{\bi{2k}k^3}{4096^k}=\f{16}{\pi},$$
by part (i) of Conjecture \ref{abm} we should have
 \begin{equation*}\label{4096Hk}
 \sum_{k=0}^\infty\f{\bi{2k}k^3}{4096^k}\l((42k+5)(H_{2k}-H_k)+7\r)=\f{32\log2}{\pi},
 \end{equation*}
 and
 \begin{align*}&\sum_{k=0}^{p-1}\f{\bi{2k}k^3}{4096^k}\l(6(42k+5)(H_{2k}-H_k)+42\r)
\\\eq\ &\l(\f{-1}p\r)\l(42+5p\,q_p(2^{12})\r)\eq\l(\f{-1}p\r)\l(42+60p\,q_p(2)\r)
\pmod{p^2}
\end{align*}
for any odd prime $p$.
\end{remark}

 \begin{conjecture} [{\rm 2022-10-11}] {\rm (i)} We have the identity
 \begin{equation}
 \sum_{k=0}^\infty(42k+5)\f{\bi{2k}k^3}{4096^k}\l(H_{2k}^{(2)}-\f{25}{92}H_{k}^{(2)}\r)=\f{2\pi}{69}.
 \end{equation}

 {\rm (ii)} Let $p>3$ be a prime. If $p\not=23$, then
 \begin{equation}\begin{aligned}
 &\sum_{k=0}^{(p-1)/2}(42k+5)\f{\bi{2k}k^3}{4096^k}\l(H_{2k}^{(2)}-\f{25}{92}H_{k}^{(2)}\r)
 \\\eq\ &\l(\f{-1}p\r)\f3{20}\l(p^4B_{p-5}-5H_{p-1}\r)\pmod{p^5}.
 \end{aligned}\end{equation}
 Also,
  \begin{equation}
 \sum_{k=0}^{p-1}(42k+5)\f{\bi{2k}k^3}{4096^k}\l(H_{2k}^{(2)}-\f{25}{92}H_{k}^{(2)}\r)
 \eq-pE_{p-3}
 \pmod{p^2}.
 \end{equation}
 \end{conjecture}
 \begin{remark} It is interesting to compare this conjecture with Remark \ref{4096R}.
 \end{remark}

 \begin{conjecture} [{\rm 2022-10-12}] \label{Conj4096} {\rm (i)} We have
 \begin{equation} \sum_{k=0}^\infty(42k+5)\f{\bi{2k}k^3}{4096^k}\l(H_{2k}^{(3)}-\f{43}{352}H_k^{(3)}\r)
 =\f{555}{77}\cdot\f{\zeta(3)}{\pi}-\f{32}{11}G.
 \end{equation}

 {\rm (ii)} For any prime $p>7$, we have
 \begin{equation}\sum_{k=0}^{(p-1)/2}(42k+5)\f{\bi{2k}k^3}{4096^k}
 \l(11H_{2k}^{(3)}-\f{43}{32}H_k^{(3)}\r)\eq-16E_{p-3}\pmod p
 \end{equation}
 and
 \begin{equation}\sum_{k=0}^{p-1}(42k+5)\f{\bi{2k}k^3}{4096^k}
 \l(11H_{2k}^{(3)}-\f{43}{32}H_k^{(3)}\r)\eq-27E_{p-3}\pmod p.
 \end{equation}
 \end{conjecture}
 \begin{remark} Conjecture \ref{Conj4096} looks quite challenging.
 \end{remark}

 \begin{conjecture} [{\rm 2022-12-05}] {\rm (i)} We have
  \begin{equation}\label{216'Hk}
 \sum_{k=0}^\infty\f{\bi{2k}k^2\bi{3k}k}{216^k}\l((6k+1)(H_{2k}-2H_k)+3\r)=\f{9\sqrt3\log3}{2\pi}
 \end{equation}
 and
 \begin{equation}\label{216Hk}
 \sum_{k=0}^\infty\f{\bi{2k}k^2\bi{3k}k}{216^k}(6k+1)(3H_{3k}-H_k)=\f{9\sqrt3\log2}{\pi}.
 \end{equation}

 {\rm (ii)} For any prime $p>3$, we have
 \begin{equation}\sum_{k=0}^{p-1}\f{\bi{2k}k^2\bi{3k}k}{216^k}((6k+1)(H_{2k}-2H_k)+3)
 \eq\l(\f p3\r)\f{3^p+3}2\pmod{p^2}
 \end{equation}
 and
 \begin{equation}\sum_{k=0}^{p-1}\f{\bi{2k}k^2\bi{3k}k}{216^k}(6k+1)(3H_{3k}-H_k)
 \eq3\l(\f p3\r)p\,q_p(2)\pmod{p^2}.
 \end{equation}
 \end{conjecture}
 \begin{remark} This is motivated by the Ramanujan series (cf. \cite[Chapter 14]{Co17} and \cite{R})
 $$\sum_{k=0}^\infty(6k+1)\f{\bi{2k}k^2\bi{3k}k}{216^k}=\f{3\sqrt3}{\pi}.$$
 \end{remark}

  \begin{conjecture} [{\rm 2022-12-04}]\label{-1024} {\rm (i)} We have
 \begin{equation}\label{1024Hk}
 \sum_{k=0}^\infty\f{\bi{2k}k^2\bi{4k}{2k}}{(-1024)^k}\l((20k+3)(H_{2k}-3H_k)+12\r)=\f{56\log2}{\pi}.
 \end{equation}

 {\rm (ii)} For any odd prime $p$, we have
 \begin{equation}\sum_{k=0}^{p-1}\f{\bi{2k}k^2\bi{4k}{2k}}{(-1024)^k}
 \l((20k+3)(H_{2k}-3H_k)+12\r)
 \eq\l(\f{-1}p\r)(12+21p\,q_p(2))\pmod{p^2}.
 \end{equation}
 \end{conjecture}
 \begin{remark} This is motivated by the Ramanujan series (cf. \cite[Chapter 14]{Co17} and \cite{R})
 $$\sum_{k=0}^\infty(20k+3)\f{\bi{2k}k^2\bi{4k}{2k}}{(-1024)^k}=\f{8}{\pi}.$$
 \end{remark}

 \begin{conjecture} [{\rm 2023-06-16}] \label{1024H2} {\rm (i)} We have
 \begin{equation}\label{4k+1}\sum_{k=0}^\infty\f{\bi{2k}k^2\bi{4k}{2k}}{(-1024)^k}\l((20k+3)\l(H_{4k}^{(2)}-\f{H_{2k}^{(2)}}4-\f{H_k^{(2)}}{16}\r)
 +\f1{4k+1}\r)=\f{\pi}6.
 \end{equation}

 {\rm (ii)} For any prime $p>3$, we have
 \begin{equation}\begin{aligned}&\ p\sum_{k=0}^{p-1}\f{\bi{2k}k^2\bi{4k}{2k}}{(-1024)^k}
 \l((20k+3)\l(H_{4k}^{(2)}-\f{H_{2k}^{(2)}}4-\f{H_k^{(2)}}{16}\r)+\f1{4k+1}\r)
 \\&\qquad\eq\l(\f{-1}p\r)+\f{10}3p^2E_{p-3}\pmod{p^3}.
 \end{aligned}
 \end{equation}
 \end{conjecture}
 \begin{remark} The identity \eqref{4k+1} looks curious and challenging.
 \end{remark}

 \begin{conjecture} [{\rm 2022-12-04}] {\rm (i)} We have
 \begin{equation}\label{48^2Hk}
 \sum_{k=0}^\infty\f{\bi{2k}k^2\bi{4k}{2k}}{48^{2k}}\l((8k+1)(3H_{2k}-4H_k)+6\r)=\f{16\sqrt3\log2}{\pi}.
 \end{equation}

 {\rm (ii)} For any prime $p>3$, we have
 \begin{equation}\sum_{k=0}^{p-1}\f{\bi{2k}k^2\bi{4k}{2k}}{48^{2k}}((8k+1)(3H_{2k}-4H_k)+6)
 \eq\l(\f p3\r)(6+8p\,q_p(2))\pmod{p^2}.
 \end{equation}
 \end{conjecture}
 \begin{remark}  This is motivated by the Ramanujan series \eqref{48^2R}.
 \end{remark}

\begin{conjecture} [{\rm 2022-10-15}] For any prime $p>3$, we have
 \begin{align}
 \label{48}\sum_{k=0}^{(p-1)/2}(8k+1)\f{\bi{2k}k^2\bi{4k}{2k}}{48^{2k}}\l(H_{2k}^{(2)}-\f5{18}H_k^{(2)}\r)
 &\eq \f p{36}B_{p-2}\l(\f13\r)\pmod{p^2},
 \\\label{48+}\sum_{k=0}^{p-1}(8k+1)\f{\bi{2k}k^2\bi{4k}{2k}}{48^{2k}}\l(H_{2k}^{(2)}-\f5{18}H_k^{(2)}\r)
 &\eq -\f 5{24}pB_{p-2}\l(\f13\r)\pmod{p^2}.
 \end{align}
 \end{conjecture}
 \begin{remark}
 The congruence \eqref{48} is motivated by \eqref{48^2}.
 \end{remark}

\begin{conjecture} [{\rm 2022-10-19}] {\rm (i)} For any prime $p>3$, we have
 \begin{equation}\begin{aligned}
 &\sum_{k=0}^{(p-1)/2}(4k+1)\f{\bi{2k}k^4}{256^k}\l(H_{2k}^{(2)}-\f12H_k^{(2)}\r)
 \\\eq\ &\sum_{k=0}^{p-1}(4k+1)\f{\bi{2k}k^4}{256^k}\l(H_{2k}^{(2)}-\f12H_k^{(2)}\r)
 \eq\f 76 p^2B_{p-3}\pmod{p^3}.
 \end{aligned}\end{equation}

 {\rm (ii)} For any odd prime $p$, we have
 \begin{equation}
 \sum_{k=0}^{(p-1)/2}(4k+1)\f{\bi{2k}k^4}{256^k}H_{2k}^{(3)}
 \eq \sum_{k=0}^{p-1}(4k+1)\f{\bi{2k}k^4}{256^k}H_{2k}^{(3)}\eq\f 32pB_{p-3}\pmod{p^2}.
 \end{equation}
 \end{conjecture}
 \begin{remark} Guo and Lian \cite[(1.7)]{GL} proved that for any prime $p>3$
 we have
 $$\sum_{k=0}^{(p-1)/2}(4k+1)\f{\bi{2k}k^4}{256^k}\l(H_{2k}^{(2)}-\f12H_k^{(2)}\r)
 \eq0\pmod{p^2}.$$
 \end{remark}

\section{Series and congruences with summands
\\ containing at least five binomial coefficients}
 \setcounter{theorem}{0}
 \setcounter{proposition}{0}

 The following two conjectures are motivated by the identity
$$\sum_{k=1}^\infty\frac{(-1)^k(205k^2-160k+32)}{k^5\binom{2k}k^5}=-2\zeta(3)$$
established by T. Amdeberhan and D. Zeilberger \cite{AZ} in 1997 via the WZ method.

\begin{conjecture} [{\rm 2022-12-09}] {\rm (i)} We have
\begin{equation}\sum_{k=1}^\infty\f{(-1)^{k-1}}{k^5\bi{2k}k^5}
\l((205k^2-160k+32)(H_{2k-1}-H_{k-1})-41k+16\r)=\f{\pi^4}{60}.
\end{equation}

{\rm (ii)} For any prime $p>5$, we have
\begin{equation}\begin{aligned}&\sum_{k=0}^{p-1}(-1)^k\bi{2k}k^5
\l((205k^2+160k+32)(H_{2k}-H_{k})+41k+16\r)
\\&\qquad\qquad\eq16p+64p^2H_{p-1}\pmod{p^6}.
\end{aligned}
\end{equation}
\end{conjecture}

\begin{conjecture} [{\rm 2022-12-09}] {\rm (i)} We have
\begin{equation}\sum_{k=1}^\infty\frac{(-1)^k((205k^2-160k+32)(4H_{2k-1}^{(2)}-12H_{k-1}^{(2)})-43)}
{k^5\binom{2k}k^5}=-8\zeta(5).
\end{equation}

{\rm (ii)} Let $p>3$ be a prime. Then
\begin{equation}\begin{aligned}&\sum_{k=1}^{(p-1)/2}\f{(-1)^k((205k^2-160k+32)(4H_{2k-1}^{(2)}-12H_{k-1}^{(2)})-43)}{k^5\bi{2k}k^5}
\\&\qquad\qquad\eq-200B_{p-5}\pmod p,
\end{aligned}
\end{equation}
and
\begin{equation}
\begin{aligned}&\sum_{k=1}^{p-1}(-1)^k\bi{2k}k^5((205k^2+160k+32)(4H_{2k}^{(2)}-12H_{k}^{(2)})+43)
\\&\qquad\qquad\eq256pH_{p-1}\pmod{p^5}
\end{aligned}
\end{equation}
if $p>5$.
\end{conjecture}

 The following two conjectures are motivated by the identity
$$\sum_{k=1}^\infty\frac{(10k^2-6k+1)(-256)^k}{k^5\binom{2k}k^5}=-28\zeta(3)$$
(cf. \cite[Identity 8]{G}).

\begin{conjecture} [{\rm 2022-12-09}] {\rm (i)} We have
\begin{equation}\sum_{k=1}^\infty\f{(-256)^k}{k^5\bi{2k}k^5}
\l((10k^2-6k+1)(2H_{2k-1}-H_{k-1})-3k+1\r)=-\f{\pi^4}{2}.
\end{equation}

{\rm (ii)} For any prime $p>3$, we have
\begin{equation}\begin{aligned}&\sum_{k=0}^{(p-1)/2}\f{\bi{2k}k^5}{(-256)^k}
\l((10k^2+6k+1)(2H_{2k}-H_{k})+3k+1\r)
\\&\qquad\qquad\eq p+\f{14}3p^4B_{p-3}\pmod{p^5}.
\end{aligned}
\end{equation}
\end{conjecture}
\begin{remark} For any prime $p>3$, the author \cite[Conjecture 31(ii)]{S19} conjectured the congruence
$$\sum_{k=(p+1)/2}^{p-1}\f{\bi{2k}k^5}{(-256)^k}(10k^2+6k+1)\eq-\f 72 p^5B_{p-3}\pmod {p^6}$$
which implies that
\begin{align*}&\ \sum_{k=(p+1)/2}^{p-1}\f{\bi{2k}k^5}{(-256)^k}\l((10k^2+6k+1)(2H_{2k}-H_k)+3k+1\r)
\\\eq&\ \sum_{k=(p+1)/2}^{p-1}\f{\bi{2k}k^5}{(-256)^k}(10k^2+6k+1)\cdot\f2p\eq-7p^4B_{p-3}\pmod{p^5}.
\end{align*}
\end{remark}

\begin{conjecture} [{\rm 2022-12-09}] {\rm (i)} We have
\begin{equation}\sum_{k=1}^\infty\frac{(-256)^k((10k^2-6k+1)(4H_{2k-1}^{(2)}-3H_{k-1}^{(2)})-2)}
{k^5\binom{2k}k^5}=-124\zeta(5).
\end{equation}

{\rm (ii)} For any prime $p>3$, we have
\begin{equation}\begin{aligned}&\sum_{k=1}^{(p-1)/2}\f{(-256)^k}{k^5\bi{2k}k^5}
((10k^2-6k+1)(4H_{2k-1}^{(2)}-3H_{k-1}^{(2)})-2)
\eq-124B_{p-5}\pmod p
\end{aligned}
\end{equation}
and
\begin{equation}\begin{aligned}&\sum_{k=1}^{(p-1)/2}\f{\bi{2k}k^5}{(-256)^k}
((10k^2+6k+1)(4H_{2k}^{(2)}-3H_{k}^{(2)})+2)
\eq\f{28}3p^3B_{p-3}\pmod {p^4}.
\end{aligned}
\end{equation}
\end{conjecture}

The following two conjectures are motivated by the identity
$$\sum_{k=1}^\infty\frac{(28k^2-18k+3)(-64)^k}{k^5\binom{2k}k^4\binom{3k}k}=-14\zeta(3)$$
conjectured by the author \cite{EJC} and confirmed recently by Au \cite{Au22}.

\begin{conjecture} [{\rm 2022-12-09}] {\rm (i)} We have
\begin{equation}\sum_{k=1}^\infty\f{(-64)^k}{k^5\bi{2k}k^4\bi{3k}k}
\l((28k^2-18k+3)(4H_{2k-1}-3H_{k-1})-20k+6\r)=-\f{\pi^4}{2}.
\end{equation}

{\rm (ii)} For any odd prime $p$, we have
\begin{equation}\begin{aligned}&\sum_{k=0}^{p-1}\f{\bi{2k}k^4\bi{3k}k}{(-64)^k}
\l((28k^2+18k+3)(4H_{2k}-3H_{k})+20k+6\r)
\\&\qquad\qquad\eq 6p-14p^4B_{p-3}\pmod{p^5}.
\end{aligned}
\end{equation}
\end{conjecture}

\begin{conjecture} [{\rm 2022-12-09}] {\rm (i)} We have
\begin{equation}\sum_{k=1}^\infty\frac{(-64)^k((28k^2-18k+3)(2H_{2k-1}^{(2)}-3H_{k-1}^{(2)})-2)}
{k^5\binom{2k}k^4\binom{3k}k}=-31\zeta(5).
\end{equation}

{\rm (ii)} For any odd prime $p$, we have
\begin{equation}\begin{aligned}&\sum_{k=1}^{p-1}\f{\bi{2k}k^4\bi{3k}k}{(-64)^k}
\l((28k^2+18k+3)(2H_{2k}^{(2)}-3H_k^{(2)})+2\r)
\\&\qquad\qquad\eq-7p^3B_{p-3}\pmod{p^4}.
\end{aligned}\end{equation}
\end{conjecture}

 The following two conjectures are motivated by the identity
$$\sum_{k=0}^\infty(20k^2+8k+1)\frac{\binom{2k}k^5}{(-4096)^k}=\frac 8{\pi^2}$$
(cf. \cite[Identity 8]{G}).

\begin{conjecture} [{\rm 2022-12-09}] {\rm (i)} We have
\begin{equation}\sum_{k=0}^\infty\f{\bi{2k}k^5}{(-4096)^k}((20k^2+8k+1)H_k-6k-1)=-\f{16\log2}{\pi^2}
\end{equation}
and
\begin{equation}\sum_{k=0}^\infty\f{\bi{2k}k^5}{(-4096)^k}(5(20k^2+8k+1)H_{2k}-10k-1)=-\f{32\log2}{\pi^2}.
\end{equation}

{\rm (ii)} For any prime $p>3$, we have
\begin{equation}\begin{aligned}&\sum_{k=0}^{(p-1)/2}\f{\bi{2k}k^5}{(-4096)^k}
((20k^2+8k+1)H_k-6k-1)
\\\eq\ &-p-2p^2q_p(2)+p^3q_p(2)^2-\f 23p^4q_p(2)^3\pmod{p^5},
\end{aligned}
\end{equation}
and
\begin{equation}\begin{aligned}&\sum_{k=0}^{(p-1)/2}\f{\bi{2k}k^5}{(-4096)^k}
(5(20k^2+8k+1)H_k-10k-1)
\\\eq\ &-p-4p^2q_p(2)+2p^3q_p(2)^2\pmod{p^4}.
\end{aligned}
\end{equation}
\end{conjecture}

\begin{conjecture} [{\rm 2022-12-09}] {\rm (i)} We have
\begin{equation}\sum_{k=1}^\infty\frac{\binom{2k}k^5}{(-4096)^k}((20k^2+8k+1)(8H_{2k}^{(2)}-3H_k^{(2)})+4)=-\frac 43.
\end{equation}

{\rm (ii)} For any prime $p>3$, we have
\begin{equation}\sum_{k=1}^{(p-1)/2}\f{\bi{2k}k^5}{(-4096)^k}
((20k^2+8k+1)(8H_{2k}^{(2)}-3H_k^{(2)})+4)\eq \f{14}3p^3B_{p-3}\pmod{p^4}
\end{equation}
and
\begin{equation}\sum_{k=1}^{p-1}\f{\bi{2k}k^5}{(-4096)^k}
((20k^2+8k+1)(8H_{2k}^{(2)}-3H_k^{(2)})+4)\eq -\f{28}3p^3B_{p-3}\pmod{p^4}.
\end{equation}
\end{conjecture}

The following two conjectures are motivated by the identity
$$\sum_{k=0}^\infty(820k^2+180k+13)\frac{\binom{2k}k^5}{(-2^{20})^k}=\frac{128}{\pi^2}$$
(cf. \cite[Identity 9]{CG}).

\begin{conjecture} [{\rm 2022-12-09}] {\rm (i)} We have
\begin{equation}\sum_{k=0}^\infty\f{\bi{2k}k^5}{(-2^{20})^k}\l((820k^2+180k+13)(H_{2k}-H_k)+164k+18\r)
=\f{256\log2}{\pi^2}.
\end{equation}

{\rm (ii)} For any odd prime $p$, we have
\begin{equation}\begin{aligned}&\sum_{k=0}^{(p-1)/2}\f{\bi{2k}k^5}{(-2^{20})^k}\l((820k^2+180k+13)(H_{2k}-H_k)+164k+18\r)
\\&\quad\eq\ 18p+26p^2q_p(2)-13p^3q_p(2)^2\pmod{p^4}.
\end{aligned}
\end{equation}
\end{conjecture}

\begin{conjecture} [{\rm 2022-12-09}] {\rm (i)} We have
\begin{equation}\sum_{k=1}^\infty\frac{\binom{2k}k^5}{(-2^{20})^k}((820k^2+180k+13)(11H_{2k}^{(2)}-3H_k^{(2)})+43)=-\frac{1}3.
\end{equation}

{\rm (ii)} Let $p>3$ be a prime. Then
\begin{equation}\begin{aligned}&\sum_{k=1}^{p-1}\frac{\binom{2k}k^5}{(-2^{20})^k}((820k^2+180k+13)(11H_{2k}^{(2)}-3H_k^{(2)})+43)
\\&\qquad\qquad\eq-\f{77}6p^3B_{p-3}\pmod{p^4},
\end{aligned}\end{equation}
and
\begin{equation}\begin{aligned}&\sum_{k=1}^{(p-1)/2}\frac{\binom{2k}k^5}{(-2^{20})^k}((820k^2+180k+13)(11H_{2k}^{(2)}-3H_k^{(2)})+43)
\\&\qquad\qquad\eq-\f{11}4pH_{p-1}\pmod{p^5}
\end{aligned}
\end{equation}
if $p>5$.
\end{conjecture}

The following two conjectures are motivated by the known identity
$$\sum_{k=0}^\infty(74k^2+27k+3)\frac{\binom{2k}k^4\binom{3k}k}{4096^k}=\frac{48}{\pi^2}$$
(cf. \cite{G11}).

\begin{conjecture} [{\rm 2022-12-09}] {\rm (i)} We have
\begin{equation}
\sum_{k=0}^\infty\f{\bi{2k}k^4\bi{3k}k}{4096^k}\l((74k^2+27k+3)H_{2k}-17k-3\r)=0
\end{equation}
and
\begin{equation}\begin{aligned}
&\sum_{k=0}^\infty\frac{\binom{2k}k^4\binom{3k}k}{4096^k}((74k^2+27k+3)(51H_{3k}+250H_{2k}-153H_k)+15)
\\&\qquad\qquad=\frac{9792\log2}{\pi^2}.
\end{aligned}\end{equation}

{\rm (ii)} For any odd prime $p$, we have
\begin{equation}
\sum_{k=0}^{p-1}\f{\bi{2k}k^4\bi{3k}k}{4096^k}((74k^2+27k+3)H_{2k}-17k-3)
\eq-3p+7p^4B_{p-3}\pmod{p^5},
\end{equation}
and
\begin{equation}\begin{aligned}
&\sum_{k=0}^{p-1}\f{\bi{2k}k^4\bi{3k}k}{4096^k}((74k^2+27k+3)(51H_{3k}+250H_{2k}-153H_k)+15)
\\&\qquad \eq\ 15p+612p^2q_p(2)-306p^3q_p(2)^2\pmod{p^4}.
\end{aligned}
\end{equation}
\end{conjecture}

\begin{conjecture} [{\rm 2022-12-09}] {\rm (i)} We have
\begin{equation}
\sum_{k=1}^\infty\frac{\binom{2k}k^4\binom{3k}k}{4096^k}((74k^2+27k+3)(92H_{2k}^{(2)}-33H_k^{(2)})+112)=160.
\end{equation}

{\rm (ii)} For any odd prime $p$, we have
\begin{equation}\begin{aligned}&\sum_{k=1}^{p-1}\f{\bi{2k}k^4\bi{3k}k}{4096^k}
\l((74k^2+27k+3)(92H_{2k}^{(2)}-33H_k^{(2)})+112\r)
\\&\qquad\qquad\eq644p^3B_{p-3}\pmod{p^4}.
\end{aligned}
\end{equation}
\end{conjecture}

 The following two conjectures are motivated by the identity
$$\sum_{k=0}^\infty(120k^2+34k+3)\frac{\binom{2k}k^4\binom{4k}{2k}}{2^{16k}}=\frac{32}{\pi^2}$$
(cf. \cite[Identity 10]{G}).

\begin{conjecture} [{\rm 2022-12-09}] {\rm (i)} We have
\begin{equation}
\sum_{k=0}^\infty\frac{\binom{2k}k^4\binom{4k}{2k}}{2^{16k}}(2(120k^2+34k+3)H_{4k}-16k-1)=0
\end{equation}
and
\begin{equation}
\sum_{k=0}^\infty\frac{\binom{2k}k^4\binom{4k}{2k}}{2^{16k}}((120k^2+34k+3)(H_{2k}-2H_k)+68k+9)
=\f{128\log2}{\pi^2}.
\end{equation}

{\rm (ii)} Let $p$ be an odd prime. Then
\begin{equation}
\begin{aligned}&\sum_{k=0}^{p-1}\f{\bi{2k}k^4\bi{4k}{2k}}{2^{16k}}((120k^2+34k+3)(H_{2k}-2H_k)+68k+9)
\\\eq\ &9p+12p^2q_p(2)-6p^3q_p(2)^2\pmod{p^4},
\end{aligned}
\end{equation}
and
\begin{equation}\sum_{k=0}^{p-1}\f{\bi{2k}k^4\bi{4k}{2k}}{2^{16k}}(2(120k^2+34k+3)H_{4k}-16k-1)
\eq-p+\f{77}6p^4B_{p-3}\pmod{p^5}
\end{equation}
if $p>3$.
\end{conjecture}

\begin{conjecture} [{\rm 2022-12-09}] {\rm (i)} We have
\begin{equation}\sum_{k=1}^\infty\frac{\binom{2k}k^4\binom{4k}{2k}}{2^{16k}}((120k^2+34k+3)(23H_{2k}^{(2)}
-7H_k^{(2)})+24)=\frac{16}3.
\end{equation}

{\rm (ii)} Let $p$ be an odd prime. Then
\begin{equation}\sum_{k=1}^{p-1}\frac{\binom{2k}k^4\binom{4k}{2k}}{2^{16k}}((120k^2+34k+3)(23H_{2k}^{(2)}
-7H_k^{(2)})+24)\eq\f{161}2p^3B_{p-3}\pmod{p^4},
\end{equation}
and
\begin{equation}\begin{aligned}&\sum_{k=1}^{(p-1)/2}\frac{\binom{2k}k^4\binom{4k}{2k}}{2^{16k}}((120k^2+34k+3)(23H_{2k}^{(2)}
-7H_k^{(2)})+24)
\\&\qquad\qquad\eq-23pH_{p-1}\pmod{p^5}
\end{aligned}
\end{equation}
if $p\not=5$.
\end{conjecture}

\begin{conjecture} [{\rm 2022-12-09}] \label{24^4} {\rm (i)} We have
\begin{equation}\begin{aligned}&\sum_{k=0}^\infty\f{\bi{2k}k^3\bi{3k}k\bi{4k}{2k}}{(-24^4)^k}
\l((252k^2+63k+5)(4H_{4k}+3H_{3k}-7H_k)+504k+63\r)
\\&\qquad\qquad=\f{192\log24}{\pi^2}.
\end{aligned}
\end{equation}

{\rm (ii)} For any prime $p>3$, we have
\begin{equation}
\begin{aligned}&\sum_{k=0}^{p-1}\f{\bi{2k}k^3\bi{3k}k\bi{4k}{2k}}{(-24^4)^k}
\l((252k^2+63k+5)(4H_{4k}+3H_{3k}-7H_k)+504k+63\r)
\\&\qquad\ \eq63p+5p^2q_p(24^4)-\f 52p^3q_p(24^4)^2\pmod{p^4}.
\end{aligned}
\end{equation}
\end{conjecture}
\begin{remark} Conjecture \ref{24^4} is motivated by the identity
$$\sum_{k=0}^\infty(252k^2+63k+5)\f{\bi{2k}k^3\bi{3k}k\bi{4k}{2k}}{(-24^4)^k}=\f{48}{\pi^2}
$$
(cf. \cite{CG}).
\end{remark}

\begin{conjecture} [{\rm 2023-01-16}] \label{10^6} {\rm (i)} We have
\begin{equation}\begin{aligned}&\sum_{k=0}^\infty\f{\bi{2k}k^2\bi{3k}k^2\bi{6k}{3k}}{10^{6k}}
\l(3(532k^2+126k+9)(H_{6k}-H_k)+532k+63\r)
\\&\qquad\qquad=\f{1125\log10}{4\pi^2}.
\end{aligned}
\end{equation}

{\rm (ii)} For any odd prime $p\not=5$, we have
\begin{equation}
\begin{aligned}&\sum_{k=0}^{p-1}\f{\bi{2k}k^2\bi{3k}k^2\bi{6k}{3k}}{10^{6k}}
\l(3(532k^2+126k+9)(H_{6k}-H_k)+532k+63\r)
\\&\qquad\ \eq63p+\f 92p^2q_p(10^6)-\f 94p^3q_p(10^6)^2\pmod{p^4}.
\end{aligned}
\end{equation}
\end{conjecture}
\begin{remark} Conjecture \ref{10^6} is motivated by the identity
$$\sum_{k=0}^\infty(532k^2+126k+9)\f{\bi{2k}k^2\bi{3k}k^2\bi{6k}{3k}}{10^{6k}}=\f{375}{4\pi^2}
$$
(cf. \cite{CG}).
\end{remark}

\begin{conjecture} [{\rm 2023-01-17}] \label{2^18} {\rm (i)} We have
\begin{equation}\begin{aligned}&\sum_{k=0}^\infty\f{\bi{2k}k^2\bi{3k}k^2\bi{6k}{3k}}{(-2^{18})^k}
\l(6(1930k^2+549k+45)(H_{6k}-H_k)+3860k+549\r)
\\&\qquad\qquad=\f{6912\log2}{\pi^2}.
\end{aligned}
\end{equation}

{\rm (ii)} For any odd prime $p$, we have
\begin{equation}
\begin{aligned}&\sum_{k=0}^{p-1}\f{\bi{2k}k^2\bi{3k}k^2\bi{6k}{3k}}{(-2^{18})^k}
\l(6(1930k^2+549k+45)(H_{6k}-H_k)+3860k+549\r)
\\&\qquad\ \eq549p+45p^2q_p(2^{18})-\f {45}2p^3q_p(2^{18})^2\pmod{p^4}.
\end{aligned}
\end{equation}
\end{conjecture}
\begin{remark} Conjecture \ref{2^18} is motivated by the identity
$$\sum_{k=0}^\infty(1930k^2+549k+45)\f{\bi{2k}k^2\bi{3k}k^2\bi{6k}{3k}}{(-2^{18})^k}=\f{384}{\pi^2}
$$
(cf. \cite{CG}).
\end{remark}

\begin{conjecture} [{\rm 2023-01-17}] \label{5^3} {\rm (i)} We have
\begin{equation}\begin{aligned}&\sum_{k=0}^\infty\f{\bi{2k}k^2\bi{3k}k^2\bi{6k}{3k}}{(-2^{18}3^65^3)^k}
\l(2(5418k^2+693k+29)(H_{6k}-H_k)+3612k+231\r)
\\&\qquad\qquad=\f{128\sqrt5}{\pi^2}\log(2^63^25).
\end{aligned}
\end{equation}

{\rm (ii)} For any prime $p>5$, we have
\begin{equation}
\begin{aligned}&\sum_{k=0}^{p-1}\f{\bi{2k}k^2\bi{3k}k^2\bi{6k}{3k}}{(-2^{18}3^65^3)^k}
\l(2(5418k^2+693k+29)(H_{6k}-H_k)+3612k+231\r)
\\&\ \eq\l(\f 5p\r)\l(231p+\f{29}3p^2q_p(2^{18}3^65^3)-\f{29}6p^3q_p(2^{18}3^65^3)^2\r)\pmod{p^4}.
\end{aligned}
\end{equation}
\end{conjecture}
\begin{remark} Conjecture \ref{5^3} is motivated by the identity
$$\sum_{k=0}^\infty(5418k^2+693k+29)\f{\bi{2k}k^2\bi{3k}k^2\bi{6k}{3k}}{(-2^{18}3^65^3)^k}=\f{128\sqrt5}{\pi^2}
$$
(cf. \cite{CG}).
\end{remark}

\begin{conjecture} [{\rm 2023-01-17}] \label{2^22} For $k\in\N$, set
$$H(k):=6H_{6k}+4H_{4k}-3H_{3k}-2H_{2k}-5H_k.$$

{\rm (i)} We have
\begin{equation}\begin{aligned}&\sum_{k=0}^\infty\f{\bi{2k}k^2\bi{3k}k\bi{4k}{2k}\bi{6k}{3k}}{(-2^{22}3^3)^k}
\l((1640k^2+278k+15)H(k)+3280k+278\r)
\\&\qquad\qquad=\f{256}{\sqrt3\pi^2}\log(2^{22}3^3).
\end{aligned}
\end{equation}

{\rm (ii)} For any prime $p>3$, we have
\begin{equation}
\begin{aligned}&\sum_{k=0}^{p-1}\f{\bi{2k}k^2\bi{3k}k\bi{4k}{2k}\bi{6k}{3k}}{(-2^{22}3^3)^k}
\l((1640k^2+278k+15)H(k)+3280k+278\r)
\\&\ \eq\l(\f 3p\r)\l(278p+15p^2q_p(2^{22}3^3)-\f{15}2p^3q_p(2^{22}3^3)^2\r)\pmod{p^4}.
\end{aligned}
\end{equation}
\end{conjecture}
\begin{remark} Conjecture \ref{2^22} is motivated by the identity
$$\sum_{k=0}^\infty(1640k^2+278k+15)\f{\bi{2k}k^2\bi{3k}k\bi{4k}{2k}\bi{6k}{3k}}{(-2^{22}3^3)^k}=\f{256}{\sqrt3\pi^2}
$$
(cf. \cite{CG}).
\end{remark}

\begin{conjecture} [{\rm 2023-01-17}] \label{7^4}  For $k\in\N$, set
$${\mathcal H}(k):=4H_{8k}-2H_{4k}+H_{2k}-3H_k.$$

{\rm (i)} We have
\begin{equation}\begin{aligned}&\sum_{k=0}^\infty\f{\bi{2k}k^3\bi{4k}{2k}\bi{8k}{4k}}{(2^{18}7^4)^k}
\l((1920k^2+304k+15){\mathcal H}(k)+1920k+152\r)
\\&\qquad\qquad=\f{56\sqrt7}{\pi^2}(9\log2+2\log7).
\end{aligned}
\end{equation}

{\rm (ii)} For any odd prime $p\not=7$, we have
\begin{equation}
\begin{aligned}&\sum_{k=0}^{p-1}\f{\bi{2k}k^3\bi{4k}{2k}\bi{8k}{4k}}{(2^{18}7^4)^k}
\l((1920k^2+304k+15){\mathcal H}(k)+1920k+152\r)
\\\eq\ &\l(\f 7p\r)\l(152p+\f{15}2p^2q_p(2^{18}7^4)-\f{15}4p^3q_p(2^{18}7^4)^2\r)\pmod{p^4}.
\end{aligned}
\end{equation}
\end{conjecture}
\begin{remark} Conjecture \ref{7^4} is motivated by the identity
$$\sum_{k=0}^\infty(1920k^2+304k+15)\f{\bi{2k}k^3\bi{4k}{2k}\bi{8k}{4k}}{(2^{18}7^4)^k}=\f{56\sqrt7}{\pi^2}
$$
(cf. \cite{CG}).
\end{remark}

\begin{conjecture} [{\rm 2022-12-09}] {\rm (i)} We have
\begin{equation}\sum_{k=1}^\infty\f{256^k}{k^7\bi{2k}k^7}
\l((21k^3-22k^2+8k-1)(4H_{2k-1}^{(2)}-5H_{k-1}^{(2)})-6k+2\r)=\f{\pi^6}{24}.
\end{equation}

{\rm (ii)} For any odd prime $p$, we have
\begin{equation}\begin{aligned}&\sum_{k=0}^{(p-1)/2}\f{\bi{2k}k^7}{256^k}
\l((21k^3+22k^2+8k+1)(4H_{2k}^{(2)}-5H_k^{(2)})+6k+2\r)
\\&\qquad\qquad\eq2p\pmod{p^5}.
\end{aligned}
\end{equation}
\end{conjecture}
\begin{remark} This is motivated by the identity
$$\sum_{k=1}^\infty\f{(21k^3-22k^2+8k-1)256^k}{k^7\bi{2k}k^7}=\f{\pi^4}8$$
conjectured by Guillera \cite{G03}.
\end{remark}

The following three conjectures are motivated by the identity \begin{equation}
\label{Pi^3}\sum_{k=0}^\infty(168k^3+76k^2+14k+1)\f{\bi{2k}k^7}{2^{20k}}=\f{32}{\pi^3}
\end{equation}
conjectured by B. Gourevich (cf. \cite{CG}).

\begin{conjecture} [{\rm 2022-12-09}] {\rm (i)} We have
\begin{equation}\begin{aligned}&\sum_{k=0}^\infty\f{\bi{2k}k^7}{2^{20k}}
\l(7(168k^3+76k^2+14k+1)(H_{2k}-H_k)+252k^2+76k+7\r)
\\&\qquad\qquad=\f{320\log2}{\pi^3}.
\end{aligned}
\end{equation}

{\rm (ii)} For any prime $p>5$, we have
\begin{equation}\begin{aligned}&\sum_{k=0}^{(p-1)/2}\f{\bi{2k}k^7}{2^{20k}}
\l(7(168k^3+76k^2+14k+1)(H_{2k}-H_k)+252k^2+76k+7\r)
\\&\ \eq\ \l(\f{-1}p\r)\l(7p^2+10p^3q_p(2)-5p^4q_p(2)^2+\f{10}3p^5q_p(2)^3-\f 52 p^6q_p(2)^4\r)
\\&\qquad\qquad\pmod{p^7}.
\end{aligned}
\end{equation}
\end{conjecture}

\begin{conjecture} [{\rm 2022-12-09}] {\rm (i)} We have
\begin{equation}\sum_{k=0}^\infty\f{\bi{2k}k^7}{2^{20k}}\l((168k^3+76k^2+14k+1)(16H_{2k}^{(2)}-5H_k^{(2)})
+8(6k+1)\r)=\f{80}{3\pi}.
\end{equation}

{\rm (ii)} For any prime $p>5$, we have
\begin{equation}\begin{aligned}&\sum_{k=0}^{(p-1)/2}\f{\bi{2k}k^7}{2^{20k}}\l(
(168k^3+76k^2+14k+1)(16H_{2k}^{(2)}-5H_k^{(2)})+8(6k+1)\r)
\\&\qquad\qquad\eq \l(\f{-1}p\r)8p\pmod{p^6}.
\end{aligned}\end{equation}
\end{conjecture}

\begin{conjecture} [{\rm 2023-06-19}] \label{2^20} {\rm (i)} We have
 \begin{equation}\sum_{k=0}^\infty\f{\bi{2k}k^7}{2^{20k}}\l((168k^3+76k^2+14k+1)(128H_{2k}^{(4)}-7H_k^{(4)})+\f {64}{2k+1}\r)
 =\f{976}{45}\pi.
 \end{equation}

 {\rm (ii)} For any odd prime $p\not=5$, we have
 \begin{equation}\begin{aligned}&\sum_{k=0}^{(p-3)/2}\f{\bi{2k}k^7}{2^{20k}}\l((168k^3+76k^2+14k+1)\l(128H_{2k}^{(4)}-7H_k^{(4)}\r)+\f{64}{2k+1}\r)
 \\&\qquad\qquad\eq -256 q_p(2)\pmod p.
 \end{aligned}
 \end{equation}
 \end{conjecture}

\begin{conjecture} [{\rm 2023-06-19}] \label{-2^24} Set
$$P(x)=4528x^4+3180x^3+972x^2+147x+9.$$

{\rm (i)} We have
\begin{equation}\sum_{k=0}^\infty\f{\bi{2k}k^7\bi{3k}k\bi{4k}{2k}}{(-2^{24})^k}\l(\l(H_{4k}+\f{H_{2k}}2\r)P(k)
-484k^3-108k^2+9k+3\r)=0,
\end{equation}
\begin{equation}\sum_{k=0}^\infty\f{\bi{2k}k^7\bi{3k}k\bi{4k}{2k}}{(-2^{24})^k}\l(\l(11H_{2k}^{(2)}-4H_{k}^{(2)}\r)P(k)
+1780k^2+633k+63\r)=\f{512}{\pi^2},
\end{equation}
\begin{equation}\sum_{k=0}^\infty\f{\bi{2k}k^7\bi{3k}k\bi{4k}{2k}}{(-2^{24})^k}\l(\l(9H_{2k}^{(3)}-2H_{k}^{(3)}\r)P(k)
-\f{916k^2+65k-9}{2k+1}\r)=0,
\end{equation}
and
\begin{equation}\sum_{k=0}^\infty\f{\bi{2k}k^7\bi{3k}k\bi{4k}{2k}}{(-2^{24})^k}\l(\l(49H_{2k}^{(4)}-2H_{k}^{(4)}\r)P(k)
+\f{5996k^2+3071k+329}{(2k+1)^2}\r)=\f{4096}{15}.
\end{equation}

{\rm (ii)} For any prime $p>3$, we have
\begin{equation}\begin{aligned}&\sum_{k=0}^{p-1}\f{\bi{2k}k^7\bi{3k}k\bi{4k}{2k}}{(-2^{24})^k}\l(\l(H_{4k}+\f{H_{2k}}2\r)P(k)
-484k^3-108k^2+9k+3\r)
\\&\qquad\qquad\eq 3p^3-\f{1953}{10}p^8B_{p-5}\pmod{p^9},
\end{aligned}
\end{equation}
\begin{equation}\begin{aligned}&\sum_{k=0}^{p-1}\f{\bi{2k}k^7\bi{3k}k\bi{4k}{2k}}{(-2^{24})^k}\l(\l(11H_{2k}^{(2)}-4H_{k}^{(2)}\r)P(k)
+1780k^2+633k+63\r)
\\&\qquad\qquad\eq 63p^2-\f{9207}{10}p^7B_{p-5}\pmod{p^8},
\end{aligned}
\end{equation}
\begin{equation}\begin{aligned}&\sum_{k=0}^{(p-3)/2}\f{\bi{2k}k^7\bi{3k}k\bi{4k}{2k}}{(-2^{24})^k}\l(\l(9H_{2k}^{(3)}-2H_{k}^{(3)}\r)P(k)
-\f{916k^2+65k-9}{2k+1}\r)
\\&\quad\eq9p-\f{81}2p^2H_{p-1}-\f{243}8p^6B_{p-5}\pmod{p^7},
\end{aligned}
\end{equation}
and
\begin{equation}\begin{aligned}&\sum_{k=1}^{(p-1)/2}\f{\bi{2k}k^7\bi{3k}k\bi{4k}{2k}}{(-2^{24})^k}
\l(\l(49H_{2k}^{(4)}-2H_{k}^{(4)}\r)P(k)
+\f{5996k^2+3071k+329}{(2k+1)^2}\r)
\\&\quad\eq-\f{441}2\l(pH_{p-1}+\f{3}4p^5B_{p-5}\r)\pmod{p^6}.
\end{aligned}
\end{equation}
\end{conjecture}
\begin{remark} This is motivated by the conjectural identity
$$\sum_{k=0}^\infty P(k)\f{\bi{2k}k^7\bi{3k}k\bi{4k}{2k}}{(-2^{24})^k}=\f{768}{\pi^4}$$
(cf. \cite{CG}).
\end{remark}

\begin{conjecture} [{\rm 2023-06-19}] \label{2^32} Set
$$Q(x)=43680k^4+20632k^3+4340k^2+466k+21$$
and $$R(x)=87360x^3+30948x^2+4340x+233.$$

{\rm (i)} We have
\begin{equation}\sum_{k=0}^\infty\f{\bi{2k}k^8\bi{4k}{2k}}{2^{32k}}\l((2(H_{4k}+3H_{2k}-4H_k)Q(k)
+R(k)\r)=\f{2^{15}\log 2}{\pi^4}
\end{equation}
and
\begin{equation}\sum_{k=0}^\infty\f{\bi{2k}k^8\bi{4k}{2k}}{2^{32k}}\l(\l(7H_{2k}^{(2)}-2H_{k}^{(2)}\r)Q(k)
+3624k^2+926k+69\r)=\f{2048}{3\pi^2}.
\end{equation}

{\rm (ii)} Let $p$ be an odd prime. If $p\not=5$, then
\begin{equation}\begin{aligned}&\sum_{k=0}^{p-1}\f{\bi{2k}k^8\bi{4k}{2k}}{2^{32k}}\l((2(H_{4k}+3H_{2k}-4H_k)Q(k)
+R(k)\r)
\\&\eq 233p^3+336p^4q_p(2)-168p^5q_p(2)^2+112p^6q_p(2)^3-84p^7q_p(2)^4\pmod{p^8}.
\end{aligned}
\end{equation}
Provided $p>3$, we have
\begin{equation}\begin{aligned}&\sum_{k=0}^{p-1}\f{\bi{2k}k^8\bi{4k}{2k}}{2^{32k}}\l(\l(7H_{2k}^{(2)}-2H_{k}^{(2)}\r)Q(k)
+3624k^2+926k+69\r)
\\&\qquad\eq 69p^2+\f{1953}{20}p^7B_{p-5}\pmod{p^8}.
\end{aligned}
\end{equation}
\end{conjecture}
\begin{remark} This is motivated by the conjectural identity
$$\sum_{k=0}^\infty Q(k)\f{\bi{2k}k^8\bi{4k}{2k}}{2^{32k}}=\f{2048}{\pi^4}$$
(cf. \cite{CG}).
\end{remark}

\Ack. The work was supported by the Natural Science Foundation of China (grant no. 12371004).

\end{document}